\documentclass[12pt]{amsart}

\usepackage{amssymb}
\usepackage{verbatim}
\usepackage{mathrsfs}

\newtheorem{thm}{Theorem}[section]

\newtheorem{lem}[thm]{Lemma}

\theoremstyle{definition}

\theoremstyle{remark}

\newtheorem*{rem}{Remark}

\numberwithin{equation}{section}

\newcommand{\kc}{\mathbf{k}}
\newcommand{\N}{\mathbf{N}}
\newcommand{\C}{\mathbf{C}}

\newcommand{\M}{\mathcal{M}}
\newcommand{\ET}{\mathcal{ET}}
\newcommand{\Mod}[1]{\ (\textup{mod}\ #1)}

\providecommand{\sgn}{\operatorname{sgn}}

\providecommand{\sgn}{\operatorname{sgn}}
\providecommand{\sym}{\operatorname{sym}}

\DeclareMathOperator{\res}{res}

\def \hyp {{}_2F_1}

\begin{document}

\title[Mixed moment of $GL(2)$ and $GL(3)$ $L$-functions]{Mixed moment of $GL(2)$ and $GL(3)$ $L$-functions}

\author[O. Balkanova]{Olga  Balkanova}
\address{Department of Mathematical Sciences, University of Gothenburg
and Chalmers University of Technology, SE-412 96 G\"{o}teborg, Sweden}
\email{olgabalkanova@gmail.com}

\author[G. Bhowmik]{Gautami Bhowmik}
\address{Laboratoire Painlev\'{e} LABEX-CEMPI, Universit\'{e} Lille 1, 59655
Villeneuve d’Ascq Cedex, France}
\email{bhowmik@math.univ-lille1.fr}

\author[D. Frolenkov]{Dmitry  Frolenkov}
\address{Steklov Mathematical Institute of Russian Academy of Sciences,
8 Gubkina st., Moscow, 119991, Russia}
\email{frolenkov@mi.ras.ru}

\author[N. Raulf]{Nicole Raulf}
\address{Laboratoire Painlev\'{e} LABEX-CEMPI, Universit\'{e} Lille 1,
59655 Villeneuve d’Ascq Cedex, France}
\email{nicole.raulf@math.univ-lille1.fr}

\begin{abstract}
Let $ \mathfrak{f} $ run over the space $ H_{4k} $ of primitive cusp forms of level one and weight $ 4k $, $ k \in \N $.
We prove an explicit formula for the mixed moment of the Hecke
$ L $-function $ L(\mathfrak{f}, 1/2) $ and the symmetric square $L$-function
$ L(\sym^2\mathfrak{f}, 1/2)$, relating it to the dual mixed moment of the double Dirichlet series and the Riemann zeta function weighted by the ${}_3F_{2}$ hypergeometric function. Analysing the corresponding special functions by the means of the Liouville-Green approximation followed by the saddle point method, we prove that the initial mixed moment is bounded by $\log^3k$.

\end{abstract}


\maketitle

\tableofcontents


\section{Introduction}
Let $H_{2k}$ be the normalised Hecke basis for the space of holomorphic cusp forms of even weight $2k \geq 2$ with respect to the full modular group.
Every function $ \mathfrak{f} \in H_{2k} $ has a Fourier expansion of the form
\begin{equation}
\mathfrak{f}(z)=\sum_{n\geq 1}\lambda_{\mathfrak{f}}(n)n^{k-1/2}\exp(2\pi inz),
\quad \lambda_{\mathfrak{f}}(1)=1.
\end{equation}

Consider the mixed moment at the critical point:
\begin{equation}\label{mixmoment}
\M(0,0):=\sum_{\mathfrak{f}\in H_{4k}} \omega(\mathfrak{f}) L(\mathfrak{f},1/2)L(\sym^2\mathfrak{f}, 1/2), \quad
\omega(\mathfrak{f}) := \frac{12\zeta(2)}{(4k-1)L(\sym^2\mathfrak{f},1)},
\end{equation}
where the corresponding $L$-functions are defined for $\Re{s}>1$ as
\begin{equation}\label{L def}
L(\mathfrak{f},s):=\sum_{n=1}^{\infty}\frac{\lambda_\mathfrak{f}(n)}{n^s}, \quad L(\sym^2\mathfrak{f},s):=\zeta(2s)\sum_{n=1}^{\infty}\frac{\lambda_\mathfrak{f}(n^2)}{n^s},
\end{equation}
and admit the analytic continuation to the whole complex plane.
Note that we consider only $k$ divisible by $4$ in \eqref{mixmoment} because otherwise $ L(\mathfrak{f},1/2)$ is identically zero. 

The mixed moment \eqref{mixmoment}  {\it with an extra smooth average over weight} was studied in \cite{BBFR} by combining an explicit formula for the first moment of symmetric square $L$-functions and an approximate functional equation for the Hecke $L$-function. This approach along with the Liouville-Green method appeared to be quite effective, producing an asymptotic formula with an arbitrary power saving error term. 

However, the same problem {\it without extra smooth averaging} is much more difficult in view of analysis of non-diagonal terms.  
For this reason, we modify the methodology of \cite{BBFR}, relying entirely on the method of analytic continuation. More precisely, we prove an explicit formula for the mixed moment \eqref{mixmoment}, which contains the diagonal main term of size $\log{k}$, the non-diagonal main term of size $k^{-1/2}$ and dual mixed moments weighted by ${}_3F_{2}$ hypergeometric functions. 
In order to state explicitly non-diagonal terms in \eqref{mixmoment}, it is required to introduce the generalised Dirichlet $L$-function
\begin{equation}\label{Lgeneralised}
\mathscr{L}_n(s) := \frac{\zeta(2s)}{\zeta(s)} \sum_{q=1}^{\infty}
\frac{1}{q^s} \left(\sum_{1\leq t \leq 2q;t^2 \equiv n \Mod{4q}}1\right),
\end{equation}
and the associated double Dirichlet series
\begin{equation*}
L^{-}_{f}(s) :=
\frac{\Gamma(3/4)}{2\sqrt{\pi}}\sum_{n<0}\frac{\mathscr{L}_{n}(1/2)}{|n|^{s+1/2}}, \quad L^{-}_{g}(s):=\frac{\Gamma(3/4)}{4\sqrt{\pi}}\sum_{n<0}\frac{\mathscr{L}_{4n}(1/2)}{|n|^{s+1/2}}.
\end{equation*}
\begin{thm}\label{main thm}
For any $\epsilon>0$ the following formula holds
\begin{equation}\label{M(0,0) exact formula1} 
\M(0,0)=2\M^D(0,0)+2\M^{ND}(0,0)+\frac{1}{2\pi i}\int_{(0)}G_{2k}(0,s)ds+O\left(\frac{k^{\epsilon}}{k}\right).
\end{equation}
Here the main diagonal term is given by
\begin{multline*}
\M^D(0,0)=
\frac{\zeta(3/2)}{2}
\Biggl(\frac{\pi}{2}-3\log{2\pi}+3\gamma+
2\frac{\zeta'(3/2)}{\zeta(3/2)}
+\frac{\Gamma'}{\Gamma}(2k-1/4)+\frac{\Gamma'}{\Gamma}(2k+1/4)\Biggr).
\end{multline*}
The main non-diagonal term is smaller in size and depends on the special value of the double Dirichlet series:
\begin{equation*}
\M^{ND}(0,0)=
\frac{2^{3/2}\pi}{\Gamma(3/4)}\frac{\Gamma(2k-1/4)}{\Gamma(2k+1/4)}L^{-}_{g}(1/4).
\end{equation*}
Finally, the last term involves the product of the Riemann zeta function with the double Dirichlet series weighted by ${}_3F_{2}$ hypergeometric function:
\begin{multline*}
G_{2k}(0,s)=
\frac{2^{5/2}}{\Gamma^2(3/4)}\Gamma(2k-1/4)\Gamma(3/4-2k)
\Gamma(1/2+s)\Gamma(1/4-s)\zeta(1/2-2s)
\\\times
\Biggl(
\left(1-2^{2s-1/2}\right)L^{-}_{f}(s)-
\frac{\left(1-2^{2s+1/2}\right)}{2^{2s}}L^{-}_{g}(s)
\Biggr){}_3F_{2}\left(2k-\frac{1}{4},\frac{3}{4}-2k,\frac{1}{4}-s;\frac{1}{2},\frac{3}{4}; 1\right).
\end{multline*}
\end{thm}

\begin{rem}
The formula \eqref{M(0,0) exact formula1} contains the error term $O\left(k^{-1+\epsilon}\right)$ with the goal to shorten the statement of the main result.
This error term can be replaced by a completely explicit expression. See \eqref{M(u,v) 1}, Lemma \ref{lem:second type}  and Lemma \ref{lem:et1at0v} for details.
\end{rem}

The first and the second terms in \eqref{M(0,0) exact formula1} can be estimated as follows:
\begin{equation}
\M^D(0,0)\ll \log{k}, \quad \M^{ND}(0,0)\ll k^{-1/2}.
\end{equation}
The analysis of the third term given by the integral of $G_{2k}(0,s)$ is the core part of this paper. Using the representation \eqref{G def} for $G_{2k}(0,s)$ in terms of the Mellin transform $\hat{g}_{2k}(u,v;s)$ studied in Lemma \ref{lem:gMelTrans}, it is required to estimate the weighted mixed moment of the double Dirichlet series and the Riemann zeta function:
\begin{equation}\label{eq:dualmoment}
\int_{-\infty}^{\infty}L^{-}_{f,g}(ir)\zeta(1/2-2ir)\hat{g}_{2k}(0,0;ir)dr.
\end{equation}

 The contribution of $|r|>3k$ is negligibly small  because in this range the function $\hat{g}_{2k}(0,0;ir)$ is of rapid decay by Lemma \ref{Lem gMellin big r}. In the remaining range, we use the representation \eqref{gMellin1} for $\hat{g}_{2k}(0,0;ir)$ as an integral of the function $\Phi_{2k}(u;x)$ defined by \eqref{defphi}. Applying the Liouville-Green approximation \eqref{eq:Phi2k} for $\Phi_{2k}(u;x)$ in terms of $Y_0$ and $J_0$ Bessel functions, we prove that
for $\kc:=4k-1$ we have (see Lemma \ref{lem:LGapprox})
\begin{multline}\label{eq:LGgooir}
\hat{g}_{2k}(0,0;ir)=-2^{3/2}\pi^{1/2}
\int_0^{\pi/2}\frac{(\tan x)^{2ir}}{(\sin(2x))^{1/2}}
Y_0(\kc x)x^{1/2}dx-\\
-2^{3/2}\pi^{1/2}
\int_0^{\pi/2}\frac{(\tan x)^{2ir}}{(\sin(2x))^{1/2}}
J_0(\kc x)x^{1/2}dx+O(k^{-3/2}).
\end{multline}
We remark that taking the absolute values to estimate the integrals in \eqref{eq:LGgooir} and using standard estimates for the Bessel functions yields
$\hat{g}_{2k}(0,0;ir)\ll k^{-1/2}$, and consequently $\M(0,0) \ll k^{1/2+\epsilon}$. In order to improve these bounds, we analyse the integrals in \eqref{eq:LGgooir} further by making the partition of unity and replacing the Bessel functions with their asymptotic formulas. Consequently, it is required to study the oscillating integral (see Lemma \ref{int:saddle})
\begin{equation}\label{int:saddle2}
\int_0^{\pi/2}\frac{\beta(x)\exp(i\kc h(x))}{(\sin(2x))^{1/2}}
dx,
\end{equation}
where $\beta=\beta(x)$ is a smooth characteristic function vanishing at the end points, and
\begin{equation*}
h(x)=-x+\frac{2r}{\kc}\log(\tan x).
\end{equation*}
A possible approach to estimate the integral \eqref{int:saddle2} is the saddle point method. However, 
as $4r\rightarrow k$ we encounter the problem of two coalescing saddle points. It is known that in this case the considered integral has different behaviour in three different ranges: $r$ is small, $r$ is near $k/4$, $r$ is large. As the standard saddle point method cannot be applied in such situation, we follow instead \cite[Section 9.2]{BlHan}, which describes the method that was originally developed by Chester, Friedman and Ursell \cite{CFU}, with some additional ideas due to Bleistein \cite{Bl}. As a result, we obtain a uniform expansion of  \eqref{int:saddle2} in terms of the Airy function (see \eqref{I Int4}), which yields the following result.
\begin{lem}\label{Lem gMellin small r}
Let $\delta$ be some fixed constant such that $0<\delta<1/4$. For $\kc^{1/2-\delta}<r\le\kc$ we have
\begin{equation}\label{gMellin small r est}
\hat{g}_{2k}(0,0;ir)\ll
\frac{1}{\kc^{5/6}}\min\left(1,\frac{\kc^{1/12}}{|\kc-4r|^{1/4}}\right)
+\frac{k^{-1/4-3\delta}+k^{-1/2}}{r}.
\end{equation}
\end{lem}
This lemma is the key ingredient for the proof of the second main theorem.
\begin{thm} \label{maincor}
The following upper bound holds
\begin{equation} \label{eq:mainformula}
\M(0,0) \ll \log^{3} k.
\end{equation}

\end{thm}

The estimate of Theorem \ref{maincor} is at the edge of current technology.  A similar result in the $q$-aspect was proved by Petrow \cite{Pet}. More precisely, \cite[Theorem 1]{Pet} refines the estimate of Conrey and Iwaniec \cite{CI} for the cibic moment of central $L$-values of level $q$ cusp forms twisted by quadratic characters of conductor $q$.  See also \cite{PY, PY2,Y} for related results. The role of the dual moment \eqref{eq:dualmoment} in this case is played by the weighted fourth moment of Dirichlet $L$-functions, as shown in \cite[Theorem 2]{Pet}. 

Another related problem is the cubic moment of cusp form $L$-functions of level one and large individual weight. In this direction, Peng \cite{Peng} proved that 
\begin{equation}\label{eq:peng}
\sum_{\mathfrak{f} \in H_{4k}}\omega(\mathfrak{f})L^3(\mathfrak{f},1/2)\ll k^{\epsilon}.
\end{equation}
The method of this paper allows replacing $k^{\epsilon}$ in \eqref{eq:peng} by $\log^5{k}$. Consequently, it is possible to obtain a slightly improved subconvexity estimate
\begin{equation}
L(\mathfrak{f},1/2)\ll k^{1/3}\log^{2}{k}.
\end{equation}
The dual moment in this case is the fourth moment of the Riemann zeta function weighted by the  ${}_3F_{2}$ hypergeometric function.

The common characteristic of all mentioned moments is that in order to derive asymptotic formulas, it is required to evaluate the dual moments very carefully taking into consideration the oscillatory behaviour of the corresponding ${}_3F_{2}$ hypergeometric functions. 
Our results yield a uniform approximation of ${}_3F_{2}$ in terms of simpler functions, which may be useful for further study of these problems. 

The main difference between the mixed moment \eqref{mixmoment} and the above mentioned moments is the appearence of the double Dirichlet series. This turns out to be a specific characteristic of symmetric square $L$-functions. Similar phenomenon in the level aspect
was discovered by Iwaniec-Michel \cite{IM} and Blomer \cite{Blo}. Furthermore, it is expected that refining the asymptotic formula of Munshi-Sengupta \cite{MS} for the mixed moment in the level aspect, it should be possible to obtain a second main term of size $q^{-1/2}$ which involves special values of  a certain double Dirichlet series. 

The paper is organised as follows.  Section \ref{sec:2} is devoted to the generalised Dirichlet $L$-functions and the associated double Dirichlet series. In Section \ref{sec:3} we recall the explicit formula for the twisted first moment of symmetric square $L$-functions. Section \ref{sec:4} is the core part of the paper, containing all required estimates for special functions. Finally, in Section \ref{sec:5} we prove Theorem \ref{main thm} and Theorem \ref{maincor}.

\section{Generalised Dirichlet $L$-functions}\label{sec:2}

As stated in Theorem \ref{main thm}, the non-diagonal terms of the mixed moment can be expressed in terms of the double Dirichlet series associated to the generalised Dirichlet $L$-function \eqref{Lgeneralised}.
In this section, we gather various results related to $\mathscr{L}_n(s)$ that are required for evaluation of the non-diagonal terms.

 According to \cite[Section 1]{Byk}, the function $ \mathscr{L}_n(s) $ considered as a function of $ s $ does not vanish only if $ n \equiv 0,1 \Mod{4} $.
The completed $L$-function
\begin{equation*}
\mathscr{L}_{n}^{*}(s) =
(\pi/|n|)^{-s/2} \Gamma(s/2+1/4-\sgn{n}/4) \mathscr{L}_n(s)
\end{equation*}
satisfies the functional equation (see \cite[Proposition 3, p. 130]{Z})
\begin{equation}
\mathscr{L}_{n}^{*}(s)=\mathscr{L}_{n}^{*}(1-s).
\end{equation}
For any $ \epsilon>0 $ we have (see \cite[Lemma 4.2]{BF1})
\begin{equation}\label{eq:subcL1}
\mathscr{L}_{n}(1/2)\ll |n|^{\theta+\epsilon},
\end{equation}
where $ \theta=1/6 $ is the best known subconvexity exponent for Dirichlet $ L $-functions obtained by Conrey and Iwaniec in \cite{CI}.
It follows from \eqref{eq:subcL1} and the Phragmen-Lindel\"{o}f principle that for any $ \epsilon>0 $ and $\Re{u}>0$ the following upper bound holds
\begin{equation}\label{eq:subcL}
\mathscr{L}_{n}(1/2+u)\ll |n|^{\max(\theta(1-2\Re{u}),0)+\epsilon}.
\end{equation}

Similarly to \cite[Section 2.3]{BBFR}, we define $f$ as the combination of the Maa{\ss}-Eisenstein series of weight
$ 1/2 $ and level $ 4 $ at the cusps $ \infty $ and $ 0 $, namely
\begin{equation}
f=f(z;s):= \zeta(4s-1) 
\left(E_{\infty}(z;s;1/2) + \frac{1+i}{4^s} E_{0}(z;s;1/2)\right).
\end{equation}
Furthermore, we define $g$ as follows 
\begin{equation}
g=g(z;s):=\frac{1}{2}\left( f\left(\frac{z}{4};s\right)+f\left(\frac{z+2}{4};s\right)\right).
\end{equation}
Using \cite[(2.26), (2.27)]{BBFR}, we associate to the functions $f$  and $g$ the double Dirichlet series
\begin{equation}\label{eq:lf}
L^{+}_{f}(s) :=
\frac{\Gamma(1/4)}{2\sqrt{\pi}}\sum_{n>0}\frac{\mathscr{L}_{n}(1/2)}{n^{s+1/2}},
\quad
L^{-}_{f}(s) :=
\frac{\Gamma(3/4)}{2\sqrt{\pi}}\sum_{n<0}\frac{\mathscr{L}_{n}(1/2)}{|n|^{s+1/2}},
\end{equation}
\begin{equation}\label{eq:lg}
L^{+}_{g}(s):=\frac{\Gamma(1/4)}{4\sqrt{\pi}}\sum_{n>0}\frac{\mathscr{L}_{4n}(1/2)}{n^{s+1/2}},
\quad
L^{-}_{g}(s):=\frac{\Gamma(3/4)}{4\sqrt{\pi}}\sum_{n<0}\frac{\mathscr{L}_{4n}(1/2)}{|n|^{s+1/2}}.
\end{equation}

\begin{thm}\label{thm:lfuncteq}
The functions $L^{\pm}_{f}(s)$ and $L^{\pm}_{g}(s)$ have a meromorphic
continuation to the whole complex plane and satisfy the functional
equations
\begin{multline}
L^{+}_{g}(s) =
\frac{-\pi^{2s+2}}{\sqrt{2}\Gamma^2(1/2+s)\sin^{2}{\pi s}} 
\left(\frac{\sin \pi(-s-1/4)}{\pi}L^{+}_{f}(-s)
- \frac{L^{-}_{f}(-s)}{\Gamma^2(3/4)} \right),
\end{multline}
\begin{multline}\label{eq:lgs}
L^{-}_{g}(s) = \frac{\pi^{2s+2}}{\sqrt{2}\Gamma^2(1/2+s)\sin^{2}{\pi s}}
\left(-\frac{\sin \pi(-s+1/4)}{\pi} L^{-}_{f}(-s)
+ \frac{L^{+}_{f}(-s)}{\Gamma^2(1/4)}\right).
\end{multline}
Furthermore, $L^{\pm}_{f}(s)$ and $L^{\pm}_{g}(s)$ are holomorphic in $\C$
except for a double pole at $s=1/2$.
\end{thm}
\begin{proof}
See \cite[Theorem 2.4]{BBFR}.
\end{proof}


\begin{thm}\label{cor Lfg coeff relation}
Assume that we have the  Laurent series
\begin{equation}\label{Laurent series Lfg}
L^{\pm}_{f,g}(s+1/2)=\frac{c_{f,g}^{\pm}(-2)}{s^2}+\frac{c_{f,g}^{\pm}(-1)}{s}+O(1).
\end{equation}
Then the coefficients $c_f^{\pm}$, $c_g^{\pm}$ satisfy the following identites:
\begin{equation}\label{eq:firstidentityc}
\frac{c^{+}_{f}(-2)}{\Gamma(1/4)}-
\frac{c^{-}_{f}(-2)}{\Gamma(3/4)}=0
\end{equation}
\begin{equation}
\frac{c^{+}_{f}(-1)}{\Gamma(1/4)}-
\frac{c^{-}_{f}(-1)}{\Gamma(3/4)}+
\frac{c^{-}_{f}(-2)\pi}{\Gamma(3/4)}=0
\end{equation}
\begin{equation}\label{relation for MT 1}
\frac{16\sqrt{\pi}}{\Gamma(3/4)}\left(
c_f^{-}(-2)(1-\sqrt{2})+
c_g^{-}(-2)(\sqrt{2}-1/2)
\right)=1.
\end{equation}

\end{thm}
\begin{proof}

To evaluate the Laurent expansion for $L^{+}_{f}( s) $ we apply \cite[(2.7), (2.9), (2.10)]{BBFR}  and \cite[(9.137.11)]{GR}, getting
\begin{multline*} 
L^{+}_{f}( s+1/2) = 
\frac{C \delta \widehat{b}_{\infty,0}}{s^2} 
\frac{(2\pi)^{s+1/2}}{2} \frac{2^{3/4}}{\Gamma(s+3/4)} 
{}_2F_{1}\left(-1/4,3/4; s+3/4; 1/2\right) \\ 
+  
\frac{1}{s} \frac{(2 \pi)^{s+1/2}}{2} \Bigg(\frac{2^{3/4}}{\Gamma(s+3/4)} 
{}_2F_{1}\left(-1/4,3/4; s+3/4; 1/2\right) C \delta (\widehat{a}_{\infty,0}
+ 2 \log{\delta} \, \widehat{b}_{\infty,0}) \\  
- \frac{2^{-1/4}}{\Gamma(s+7/4)} 
{}_2F_{1}\left(3/4,3/4;s+7/4;1/2\right) C \delta \widehat{b}_{\infty,0}\Bigg)  
+ O(1). 
\end{multline*} 

Note that in our case (see \cite[(2.26), (2.27)]{BBFR})  $C=\sqrt{2},$  $\delta=1/2$ and
\begin{equation}\label{coeff:ab}
 \alpha_{\infty,0}=\gamma-\log{4\pi}, \quad \widehat{a}_{\infty,0}=(\gamma-\log{8\pi})/2,\quad b_{\infty,0}=1/2,  \quad\widehat{b}_{\infty,0}=1/4.
\end{equation}

Applying \cite[(9.131.1)]{GR} we infer 
\begin{equation*} 
\begin{split} 
{}_2F_{1} & \left(-1/4,3/4; s+3/4; 1/2\right) 
= 
\left(\frac{1}{2}\right)^{s+1/4} 
\hyp\left(s+1, s; s+3/4; 1/2\right) \\ 
&= 
\left(\frac{1}{2}\right)^{s+1/4} 
\left(1 + 
\frac{\Gamma(s+3/4)}{\Gamma(s+1) \Gamma(s)} 
\sum_{n=0}^{\infty} 
\frac{\Gamma(s+2+n) \Gamma(s+1+n)}{(n+1)! \Gamma(s+7/4+n)} 
\left(\frac{1}{2}\right)^{n+1}\right),
\end{split} 
\end{equation*} 
and therefore,
\begin{equation*} 
\begin{split} 
\frac{1}{s^2} & 
\frac{(2\pi)^{s+1/2}}{2^{1/4} \Gamma(s+3/4)} 
{}_2F_{1}\left(-1/4,3/4; s+3/4; 1/2\right)= 
\frac{1}{s^2} 
\frac{\pi^{s+1/2}}{\Gamma(s+3/4)} \\ 
& \qquad 
+ 
\frac{1}{s} 
\frac{\pi^{s+1/2}}{2 \Gamma^2(s+1)} \sum_{n=0}^{\infty} 
\frac{\Gamma(s+2+n) \Gamma(s+1+n)}{(n+1)! \Gamma(s+7/4+n)} 
\left(\frac{1}{2}\right)^{n}. 
\end{split} 
\end{equation*} 
Since 
\begin{equation*} 
\frac{\pi^{s+1/2}}{\Gamma(s+3/4)} = 
\frac{\pi^{1/2}}{\Gamma(3/4)} + \frac{\pi^{1/2}}{\Gamma(3/4)} 
\left(\log \pi - \psi(3/4)\right) s + O\left(s^2\right)
\end{equation*} 
and 
\begin{multline*}  
\lim_{s \rightarrow 0} \frac{\pi^{s+1/2}}{2 \Gamma^2(s+1)} 
\sum_{n=0}^{\infty} 
\frac{\Gamma(s+2+n) \Gamma(s+1+n)}{(n+1)! \Gamma(s+7/4+n)} 
\left(\frac{1}{2}\right)^{n} \\ 
= 
\frac{(2\pi)^{1/2}}{2^{3/2}} \sum_{n=0}^{\infty} 
\frac{\Gamma(1+n) \Gamma(1+n)}{n! \Gamma(7/4+n)} 
\left(\frac{1}{2}\right)^{n} \\ 
= 
\frac{(2\pi)^{1/2}}{2^{3/2} \Gamma(7/4)} 
\hyp\left(1, 1; 7/4; 1/2\right) 
= 
\frac{(2\pi)^{1/2}}{2^{3/2} \Gamma(3/4)} 
\left(\psi(7/8) - \psi(3/8)\right),
\end{multline*} 
we obtain 
\begin{equation*} 
\begin{split} 
L^{+}_{f}( s+1/2) 
&= 
\frac{\pi^{1/2} C \delta \widehat{b}_{\infty,0}}{\Gamma(3/4)} 
\frac{1}{s^2} 
+ 
\frac{\pi^{1/2} C \delta}{\Gamma(3/4)} 
\Big(\left(\log \pi - \psi(3/4)\right) \widehat{b}_{\infty,0} \\ 
& \quad 
+ 
\widehat{a}_{\infty,0} + 2 \log{\delta} \, 
\widehat{b}_{\infty,0}\Big) \frac{1}{s} + O(1). 
\end{split} 
\end{equation*} 


Similarly, using \cite[(9.137.11)]{GR} we have
\begin{equation*}
\begin{split}
L^{-}_{f}(s+1/2) 
&= 
\frac{(2 \pi)^{s+1/2} C \delta \widehat{b}_{\infty,0}}{2^{11/4} 
\Gamma(s+5/4) s^2} {}_2F_{1}\left(1/4,5/4;s+5/4;1/2\right) \\ 
& \qquad 
+ 
\frac{(2 \pi)^{s+1/2} C \delta (\widehat{a}_{\infty,0}
+ 2 \log{\delta} \, \widehat{b}_{\infty,0})}{2^{11/4} \Gamma(s+5/4) s} 
{}_2F_{1}\left(1/4,5/4;s+5/4;1/2\right) \\ 
& \qquad 
+ 
\frac{(2 \pi)^{s+1/2} C \delta \widehat{b}_{\infty,0}}{2s} 
\frac{2^{1/4}}{\Gamma(s+5/4)} {}_2F_{1}\left(1/4,1/4;s+5/4;1/2\right)
+ O(1). 
\end{split}
\end{equation*}
In this case \cite[(9.131.1)]{GR}  yields 
\begin{equation*} 
\begin{split} 
& {}_2F_{1}\left(1/4,5/4;s+5/4;1/2\right) 
= \\ 
& \qquad 
\left(\frac{1}{2}\right)^{s-1/4} 
+ \left(\frac{1}{2}\right)^{s-1/4}
\frac{\Gamma(s+5/4) s}{2 \Gamma^2(s+1)} \sum_{n=0}^{\infty} 
\frac{\Gamma(s+2+n) \Gamma(s+1+n)}{(n+1)! \Gamma(s+9/4+n)} 
\left(\frac{1}{2}\right)^n 
\end{split} 
\end{equation*} 
so that 
\begin{equation*}
\begin{split}
L^{-}_{f}(s+1/2) &= 
\frac{\pi^{1/2} C \delta \widehat{b}_{\infty,0}}{\Gamma(1/4)} 
\frac{1}{s^2} 
+ 
\frac{\pi^{1/2} C \delta}{\Gamma(1/4)} 
\Bigg(\left(\log \pi - \psi(5/4)\right) \widehat{b}_{\infty,0} \\ 
& \quad 
+ 
4 \widehat{b}_{\infty,0} 
+ 
\widehat{a}_{\infty,0} + 2 \log{\delta} \, \widehat{b}_{\infty,0}\Bigg) 
\frac{1}{s} + O(1). 
\end{split} 
\end{equation*} 

Now let us look at $ L^{+}_{g}( s) $. Since 
\begin{equation*} 
\begin{split} 
L^{+}_{g}( s+1/2) &= 
\frac{(2 \pi)^{s+1/2}}{2} \frac{\delta^{2s+1}}{C} 
\Bigg(
\frac{1}{s^2} \frac{2^{3/4} b_{\infty,0}}{\Gamma(s+3/4)} 
\hyp(-1/4, 3/4; s+3/4;1/2) \\ 
& \qquad 
+ 
\frac{1}{s} 
\bigg(\frac{2^{3/4} a_{\infty,0}}{\Gamma(s+3/4)} \hyp(-1/4, 3/4; s+3/4;1/2) \\ 
& \qquad 
- \frac{2^{-1/4} b_{\infty, 0}}{\Gamma(s+7/4)} \hyp(3/4, 3/4; s+7/4; 1/2) 
\bigg) 
\Bigg) 
+ O(1) 
\end{split} 
\end{equation*} 
we obtain, using the same transformations for the hypergeometric functions 
as before, 
\begin{equation*} 
\begin{split} 
L^{+}_{g}( s+1/2) &= 
\frac{1}{s^2} 
\frac{(\pi \delta^2)^{1/2} b_{\infty,0}}{\Gamma(3/4) C} \\ 
& \quad 
+ 
\frac{1}{s} \frac{(\pi \delta^2)^{1/2}}{\Gamma(3/4) C} 
\Big(\big(\log(\pi \delta^2) - \psi(3/4)\big) b_{\infty,0} 
+ a_{\infty,0}\Big) 
+ O(1). \\ 
\end{split} 
\end{equation*} 
Lastly,  we show that
\begin{equation*} 
\begin{split} 
L^{-}_{g}( s+1/2) 
&= 
\frac{1}{s^2} \frac{\delta}{C} \frac{\pi^{1/2}}{\Gamma(1/4)} 
b_{\infty, 0} \\ 
& \quad 
+ 
\frac{1}{s} \frac{\delta}{C} \frac{\pi^{1/2}}{\Gamma(1/4)} 
\left(\left(\log(\pi \delta^2) - \psi\left(1/4\right)\right) 
b_{\infty, 0} 
+ a_{\infty, 0}\right) 
+ O(1). \\
\end{split} 
\end{equation*} 

The identity \eqref{eq:firstidentityc} follows immediately from the expansion for $ L^{+}_{f}( s) $ and $ L^{-}_{f}( s) $. 
Furthermore,
\begin{multline*} 
\frac{c_f^+(-1)}{\Gamma(1/4)} - \frac{c_f^+(-1)}{\Gamma(3/4)} 
+ \frac{c_f^-(-2) \pi}{\Gamma(3/4)} = 
\frac{\pi^{1/2} C \delta}{\Gamma(1/4) \Gamma(3/4)} 
\Bigg(\left(\log \pi - \psi(3/4)\right) \widehat{b}_{\infty,0}  
+ 
\widehat{a}_{\infty,0} \\+ 2 \log{\delta} \, 
\widehat{b}_{\infty,0} 
- 
\Big(\left(\log \pi - \psi(5/4)\right) \widehat{b}_{\infty,0}  
+ 
4 \widehat{b}_{\infty,0} 
+ 
\widehat{a}_{\infty,0} + 2 \log{\delta} \, \widehat{b}_{\infty,0}\Big) 
+ 
\pi \widehat{b}_{\infty,0}\Bigg) \\ 
= 
\frac{\pi^{1/2} C \delta \widehat{b}_{\infty,0}}{\Gamma(1/4) \Gamma(3/4)} 
\Big(- \psi(3/4) + \psi(1/4) + \pi\Bigg)
= 
0 
\end{multline*} 
by \cite[(8.365.10), p.~905]{GR}. Finally, 
\begin{equation*} 
\begin{split} 
\frac{16 \sqrt{\pi}}{\Gamma(3/4)} 
& 
\left(c_f^-(-2) \left(1 - \sqrt{2}\right) 
+ c_g^{-}(-2) \left(\sqrt{2} - 1/2\right)\right) \\ 
&= 
\frac{16 \sqrt{\pi}}{\Gamma(3/4)} 
\left(\frac{\pi^{1/2} C \delta \widehat{b}_{\infty,0}}{\Gamma(1/4)} 
\left(1 - \sqrt{2}\right) 
+ 
\frac{\pi^{1/2} \delta b_{\infty, 0}}{\Gamma(1/4) C} 
\left(\sqrt{2} - 1/2\right)\right)  \\ 
&= 
8 \sqrt{2} \delta 
\left(C \widehat{b}_{\infty,0} \left(1 - \sqrt{2}\right) 
+ 
\frac{b_{\infty, 0}}{C} \left(\sqrt{2} - 1/2\right)\right) = 
1 
\end{split} 
\end{equation*} 
since according to \eqref{coeff:ab} we have
$ C = \sqrt{2} $, $ \delta = 1/2 $, $ \widehat{b}_{\infty,0} = 1/4 $ 
and $ b_{\infty, 0} = 1/2 $. 
\end{proof}

\begin{thm}\label{thm second mom Lfg}
The following estimates hold
\begin{equation}
\int_0^T|L^{\pm}_{g}(it)|^2dt\ll T(\log T)^4,\quad
\int_0^T|L^{\pm}_{f}(it)|^2dt\ll T(\log T)^4.
\end{equation}
\end{thm}
\begin{proof}
This is a direct consequence of \cite[Theorem 5.1 (iv)]{Mu}. Note that the Fourier-Whittaker expansion of the functions $f(z;s)$ and $g(z;s)$ can be found in \cite[(2.26), (2.27)]{BBFR}.
\end{proof}


\section{Explicit formula for the twisted first moment of symmetric square $L$-functions}\label{sec:3}

For $0<x<1$ and $0\le\Re{u}<2k-3/2$ let
\begin{multline}\label{defpsi}
\Psi_k(u;x) :=
x^k \frac{\Gamma(k-1/4-u/2) \Gamma(k+1/4-u/2)}{\Gamma(2k)}\\ \times
{}_2F_{1}\left(k-\frac{1}{4}-\frac{u}{2}, k+\frac{1}{4}-\frac{u}{2}; 2k; x\right),
\end{multline}
\begin{multline}\label{defphi}
\Phi_k(u;x) :=
\frac{\Gamma(k-1/4-u/2) \Gamma(3/4-k-u/2)}{\Gamma(1/2)}\\ \times
{}_2F_{1}\left(k-\frac{1}{4}-\frac{u}{2}, \frac{3}{4}-k-\frac{u}{2}; 1/2; x\right),
\end{multline}
where $ {}_2F_{1}(a,b;c;x) $ is the Gauss hypergeometric function.
For simplicity, let us introduce the following notation
\begin{equation}\label{defpsiphi0}
\Psi_k(x) :=\Psi_k(0;x),\quad
\Phi_k(x) :=\Phi_k(0;x).
\end{equation}
For $1-2k<\Delta<1/2-\Re{u}$ let
\begin{multline}\label{eq:integralI}
I_k(u;x):=\frac{1}{2\pi i}\int_{(\Delta)}\frac{\Gamma(k-1/2+w/2)}{\Gamma(k+1/2-w/2)}\Gamma(\frac{1}{2}-u-w)\\ \times \sin\left( \pi \frac{1/2+u+w}{2}\right)x^wdw.
\end{multline}

According to \cite[(5.3)]{BF1}, for  $x>2$ we have
\begin{equation}\label{eq:integralI x>2}
I_k(u;x)=(-1)^k\frac{\cos(\pi(1/4+u/2))}{2^{1/2+u}\pi^{1/2}}x\Psi_k\left(u;\frac{4}{x^2}\right).
\end{equation}

According to  \cite[(5.5)]{BF1}, for $0<x<2$ we have
\begin{equation}\label{eq:integralI x<2}
I_k(u;x)=(-1)^k\frac{\sin(\pi(1/4+u/2))}{\pi^{1/2}}x^{1/2-u}\Phi_k\left(u;\frac{x^2}{4}\right).
\end{equation}

Now we are ready to state the explicit formula for the twisted first moment of symmetric
square $L$-functions.


\begin{lem}\label{lem:EF}
For $0\le\Re{u}<4k-3/2$ we have
\begin{multline*}
\sum_{\mathfrak{f} \in H_{4k}} \omega(\mathfrak{f}) \lambda_{\mathfrak{f}}(l) L(\sym^2\mathfrak{f}, 1/2+u)
=\\
M^{D}(u,l) \delta_{l=\Box} + M^{ND}(u,l) + ET_1(u,l) + ET_2(u,l),
\end{multline*}
where
\begin{equation}
\delta_{l=\Box} =
\begin{cases}
1 & \text{if } l \text{ is a full square,} \\
0 & \text{otherwise},
\end{cases}
\end{equation}
\begin{multline} \label{eq:MT}
M^D(u,l^2)
=
\frac{\zeta(1+2u)}{l^{1/2+u}} + \sqrt{2}(2\pi)^{3u}\cos{\pi(1/4+u/2)}
\times \\
\frac{\zeta(1-2u)}{l^{1/2-u}} \frac{\Gamma(2k-1/4-u/2) \Gamma(2k+1/4-u/2)
\Gamma(1-2u)}{\Gamma(2k+1/4+u/2) \Gamma(2k-1/4+u/2)\Gamma(1-u)},
\end{multline}
\begin{equation}\label{eq:MNDT}
M^{ND}(u,l)
=
\frac{(2\pi)^{1/2+u}}{2l^{1/4-u/2}} \frac{\Gamma(2k-1/4-u/2)}{\Gamma(2k+1/4+u/2)}
\mathscr{L}_{-4l}(1/2+u),
\end{equation}
\begin{equation}\label{eq:ET1}
ET_1(u,l) =
(2\pi)^{1/2+u}\sum_{1\leq n<2\sqrt{l}}\frac{\mathscr{L}_{n^2-4l}(1/2+u)}{n^{1/2-u}}I_{2k}\left(u;\frac{n}{l^{1/2}}\right),
\end{equation}
\begin{equation}\label{eq:ET2}
ET_2(u,l) =
(2\pi)^{1/2+u}\sum_{n>2\sqrt{l}}\frac{\mathscr{L}_{n^2-4l}(1/2+u)}{n^{1/2-u}}I_{2k}\left(u;\frac{n}{l^{1/2}}\right).
\end{equation}
\end{lem}
\begin{proof}
See \cite[(2.9),~(5.6)]{BF1}.
\end{proof}


\begin{rem}
The role of the shift $u$ is to guarantee the absolute convergence of the integral \eqref{eq:integralI}.
\end{rem}

\section{Special functions}\label{sec:4}

For a function $h(x)$, we denote its Mellin transform by
\begin{equation}\label{Mellin def}
\hat{h}(s)=\int_0^{\infty}h(x)x^{s-1}dx.
\end{equation}

Let us define for $0<x<1$
\begin{equation}\label{def f}
f_{2k}(u,v;x):=\frac{x^{1/2+v}}{(1-x)^{1/2+v}}I_{2k}\left(u;\frac{2}{(1-x)^{1/2}}\right),
\end{equation}
and  $f_{2k}(u,v;x):=0$ for $x>1$.

For $0<x<\infty$ let
\begin{equation}\label{def g}
g_{2k}(u,v;x):=\frac{x^{1/2+v}}{(1+x)^{1/2+v}}I_{2k}\left(u;\frac{2}{(1+x)^{1/2}}\right).
\end{equation}

In this section we analyse the Mellin transforms of the functions $f_{2k}(u,v;x)$ and $g_{2k}(u,v;x)$.

\subsection{Mellin transform of $f_{2k}$}
\begin{lem}
For $\Re{s}>-1/2-\Re{v}$ and $\Re{v}<2k$, $0\le\Re{u}<4k-1$, the Mellin transform of the function \eqref{def f} can be written in three different ways:
\begin{equation}\label{fMellin1}
\hat{f}_{2k}(u,v;s)=
\frac{2\cos(\pi(1/4+u/2))}{2^{1/2+u}\pi^{1/2}}
\int_0^1\frac{(1-x)^{s+v-1/2}}{x^{1+v}}
\Psi_{2k}\left(u;x\right)dx,
\end{equation}
\begin{multline}\label{fMellin2}
\hat{f}_{2k}(u,v;s)=\Gamma(1/2+s+v)
\frac{1}{2\pi i}\int_{(\Delta)}\frac{\Gamma(2k-1/2+w/2)}{\Gamma(2k+1/2-w/2)}\\ \times \Gamma(\frac{1}{2}-u-w)
\sin\left( \pi \frac{1/2+u+w}{2}\right)\frac{\Gamma(1/2-v-w/2)}{\Gamma(1+s-w/2)}2^wdw,
\end{multline}
where $1-4k<\Delta<\min(1-2\Re{v},1/2-\Re{u})$, and
\begin{multline}\label{fMellin3}
\hat{f}_{2k}(u,v;s)=
\frac{2^{1/2-u}\sin(\pi(3/4+u/2))}{\pi^{1/2}}\Gamma(1/2+s+v)\\\times
\frac{\Gamma(2k-1/4-u/2)\Gamma(2k+1/4-u/2)\Gamma(2k-v)}{\Gamma(4k)\Gamma(2k+1/2+s)}\\\times
{}_3F_{2}\left(2k-\frac{1}{4}-\frac{u}{2}, 2k+\frac{1}{4}-\frac{u}{2},2k-v; 4k, 2k+1/2+s; 1\right).
\end{multline}
\end{lem}
\begin{proof}
It follows from the definition of the Mellin transform  \eqref{Mellin def} that
\begin{equation}\label{fMellin0}
\hat{f}_{2k}(u,v;s)=
\int_0^1\frac{x^{s+v-1/2}}{(1-x)^{1/2+v}}I_{2k}\left(u;\frac{2}{(1-x)^{1/2}}\right)dx.
\end{equation}
Substituting \eqref{eq:integralI x>2} into \eqref{fMellin0} we obtain \eqref{fMellin1}.

Assuming first that $\Re{u}>0$, we substitute \eqref{eq:integralI} to \eqref{fMellin0}. For $\Re{u}>0,$ $\Re{w}<1-2\Re{v},$ $\Re{s}>-1/2-\Re{v}$, the resulting double integral converges absolutely. Changing the order of integration and using  \cite[(5.12.1)]{HMF}, namely
\begin{equation*}
\int_0^1\frac{x^{s+v-1/2}}{(1-x)^{1/2+w/2+v}}dx=
\Gamma(1/2+s+v)\frac{\Gamma(1/2-v-w/2)}{\Gamma(1+s-w/2)},
\end{equation*}
we obtain \eqref{fMellin2}. Note that the integral on the right-hand side of \eqref{fMellin2} converges absolutely provided that $\Re{s}>-1/2-\Re{u}-\Re{v}.$ 

Moving the line of integration in \eqref{fMellin2} to the left and crossing the poles at $w=1-4k-2j$, we finally prove \eqref{fMellin3}.
\end{proof}

\begin{lem}
For $-1/4<\Re{v}<2k$ we have
\begin{equation}\label{fMellin -1/4}
\hat{f}_{2k}(0,v;-1/4)=
\frac{\Gamma^2(1/4+v)}{\pi^{1/2}}
\frac{\Gamma(2k-1/4)\Gamma(2k-v)}{\Gamma(2k+1/4)\Gamma(2k+v)}.
\end{equation}
\end{lem}
\begin{proof}
Rewriting \eqref{fMellin3} for $u=0$, we obtain
\begin{multline*}
\hat{f}_{2k}(0,v;-1/4)=
\frac{\Gamma(1/4+v)}{\pi^{1/2}}
\frac{\Gamma(2k-1/4)\Gamma(2k-v)}{\Gamma(4k)} 
{}_2F_{1}\left(2k-\frac{1}{4}, 2k-v; 4k; 1\right).
\end{multline*}
Then \eqref{fMellin -1/4} follows by applying \cite[(15.4.20)]{HMF}.
\end{proof}

\begin{lem}\label{Lem fMellin 1/2}
The following estimates hold
\begin{equation}\label{fMellin 1/2}
\hat{f}_{2k}(0,0;1/2), \quad \frac{\partial}{\partial s}\hat{f}_{2k}(0,0;s)\Bigg|_{s=1/2}\quad
\ll\frac{k^{\epsilon}}{k^2}.
\end{equation}
\end{lem}
\begin{proof}
To prove \eqref{fMellin 1/2} we apply \eqref{fMellin1} together with the Liouville-Green approximation of the function $\Psi_{2k}\left(x\right)$ obtained in \cite{BF1}. More precisely, using \cite[(6.58), (6.62), (6.64), (6.68)]{BF1}, we have that
\begin{multline}\label{psik LG}
\Psi_{2k}\left(\frac{1}{\cosh^2{\sqrt{\xi}/2}} \right)\left( \xi\sinh^2{\sqrt{\xi}}\right)^{1/4}=\\
C_K\left(
\sqrt{\xi}K_0(u\sqrt{\xi})-\frac{\xi}{u}K_1(u\sqrt{\xi})B_K(0;\xi)\right)+
O\left(\frac{\sqrt{\xi}K_0(u\sqrt{\xi})}{u^{3}}\min\left(\sqrt{\xi}, \frac{1}{\xi}\right)\right),
\end{multline}
where $u=2k-1/2,$ $C_K=2+O(k^{-1})$  and
\begin{equation*}
B(0;\xi)=\frac{1}{16}\left( \frac{\coth{\sqrt{\xi/4}}}{\sqrt{\xi}}-\frac{2}{\xi}\right).
\end{equation*}
Note that there is a typo in the formula \cite[(6.58)]{BF1} for $B(0;\xi).$ Instead of $\coth{\sqrt{\xi}}$ there should be $\coth{\sqrt{\xi/4}}.$
It follows from \eqref{psik LG} and the standard bounds on the $K$-Bessel functions \cite[(10.25.3), (10.30.2), (10.30.3)]{HMF} that
\begin{equation}\label{psik LG2}
\Psi_{2k}\left(\frac{1}{\cosh^2{\sqrt{\xi}/2}} \right)\left( \xi\sinh^2{\sqrt{\xi}}\right)^{1/4}\ll
\sqrt{\xi}K_0(u\sqrt{\xi}).
\end{equation}
Applying  \eqref{fMellin1} and making the change of variable $x=\cosh^{-2}{\sqrt{\xi}/2}$, we obtain
\begin{multline*}
\hat{f}_{2k}(0,0;1/2)\ll
\int_0^1x^{-1}\Psi_{2k}\left(x\right)dx\ll
\int_0^{\infty}\Psi_{2k}\left(\frac{1}{\cosh^2{\sqrt{\xi}/2}} \right)\frac{\sinh{\sqrt{\xi}/2}}{\cosh{\sqrt{\xi}/2}}
\frac{d\xi}{\xi^{1/2}}.
\end{multline*}
Then according to \eqref{psik LG2} we have
\begin{equation}\label{fMellin 1/2 2}
\hat{f}_{2k}(0,0;1/2)\ll
\int_0^{\infty}|K_0(u\sqrt{\xi})|
\frac{\tanh{\sqrt{\xi}/2}}{\sinh^{1/2}{\sqrt{\xi}}}
\frac{d\xi}{\xi^{1/4}}.
\end{equation}
Estimating the $K$-Bessel function by the means of \cite[(10.25.3), (10.30.3)]{HMF} completes the proof of the first estimate in  \eqref{fMellin 1/2}.
The derivative of $\hat{f}_{2k}(0,0;s)$ can be estimated similarly since it follows from \eqref{fMellin1} that
\begin{equation*}
\frac{\partial}{\partial s}\hat{f}_{2k}(0,0;s)\Bigg|_{s=1/2}\ll
\int_0^1\frac{\log(1-x)}{x}\Psi_{2k}\left(x\right)dx.
\end{equation*}
\end{proof}

\begin{lem}
The following estimate holds
\begin{equation}\label{fMellin ir}
\hat{f}_{2k}(0,0;ir)\ll\frac{k^{\epsilon}}{k(1+|r|)^2}.
\end{equation}
\end{lem}
\begin{proof}
For $|r|\ll1$  we estimate \eqref{fMellin1} trivially:
\begin{equation*}
\hat{f}_{2k}(0,0;ir)\ll
\int_0^1\frac{(1-x)^{-1/2}}{x}\Psi_{2k}\left(x\right)dx.
\end{equation*}
Repeating the arguments of Lemma \ref{Lem fMellin 1/2}, we obtain
\begin{equation*}
\hat{f}_{2k}(0,0;ir)\ll
\int_0^{\infty}|K_0(u\sqrt{\xi})|
\frac{d\xi}{\xi^{1/4}\sinh^{1/2}{\sqrt{\xi}}}.
\end{equation*}
Using \cite[(10.25.3), (10.30.2)]{HMF} we prove  \eqref{fMellin ir}.

Now let us consider the case $|r|\gg 1$.
Introducing the notation
\begin{equation}\label{Y def}
T_{2k}(x):=(1-x)^{1/2}\Psi_{2k}(x),
\end{equation}
we have
\begin{equation}\label{fMellinY}
\hat{f}_{2k}(0,0;ir)=
\frac{1}{\pi^{1/2}}\int_0^1\frac{(1-x)^{ir-1}}{x}T_{2k}(x)dx.
\end{equation}
Integrating  \eqref{fMellinY} by parts three times, we obtain
\begin{equation}\label{fMellinY2}
\hat{f}_{2k}(0,0;ir)\ll\frac{1}{(1+|r|)^3}
\int_0^1(1-x)^{ir+2}\left(T_{2k}(x)x^{-1}\right)'''dx.
\end{equation}

According to  \cite[(6.47)]{BF1}, the function $T_{2k}(x)$ satisfies the differential equation 
\begin{equation}\label{eq:diffurf21}
T_{2k}''(x)-(u^2\alpha(x)+\beta(x))T_{2k}(x)=0,
\end{equation}
where $u=2k-1/2$ and
\begin{equation}\label{def fandg}
\alpha(x):=\frac{1}{x^2(1-x)}, \quad \beta(x):=-\frac{1}{4x^2(1-x)^2}+\frac{3}{16x(1-x)}.
\end{equation}
Differentiating \eqref{eq:diffurf21} yields
\begin{equation*}
T_{2k}'''(x)=(u^2\alpha'(x)+\beta'(x))T_{2k}(x)+(u^2\alpha(x)+\beta(x))T'_{2k}(x).
\end{equation*}
Consequently,
\begin{multline}\label{3rd derivative}
\left(T_{2k}(x)x^{-1}\right)'''=\left(\frac{u^2\alpha(x)+\beta(x)}{x}+\frac{6}{x^3}\right)T'_{2k}(x)\\+
\left(\frac{u^2\alpha'(x)+\beta'(x)}{x}-3\frac{u^2\alpha(x)+\beta(x)}{x^2}-\frac{6}{x^4}\right)T_{2k}(x).
\end{multline}
Substituting \eqref{3rd derivative} into \eqref{fMellinY2}, we have
\begin{multline}\label{fMellinY3}
\hat{f}_{2k}(0,0;ir)\ll\frac{1}{(1+|r|)^3}\\\times
\int_0^1(1-x)^{2}\left(\frac{u^2|\alpha'(x)|+|\beta'(x)|}{x}+\frac{u^2|\alpha(x)|+|\beta(x)|}{x^2}+\frac{1}{x^4}\right)|T_{2k}(x)|dx\\+
\frac{1}{(1+|r|)^3}\left|\int_0^1(1-x)^{2+ir}
\left(\frac{u^2\alpha(x)+\beta(x)}{x}+\frac{6}{x^3}\right)T'_{2k}(x)dx\right|.
\end{multline}
To estimate the second integral in \eqref{fMellinY3}, we integrate it by parts, getting
\begin{multline}\label{fMellinY4}
\frac{1}{(1+|r|)^3}\int_0^1(1-x)^{2+ir}
\left(\frac{u^2\alpha(x)+\beta(x)}{x}+\frac{6}{x^3}\right)T'_{2k}(x)dx\\\ll
\frac{1}{(1+|r|)^2}\int_0^1(1-x)
\left(\frac{u^2|\alpha(x)|+|\beta(x)|}{x}+\frac{1}{x^3}\right)|T_{2k}(x)|dx+\\
\int_0^1 \frac{(1-x)^2}{(1+|r|)^3}
\left(\frac{u^2|\alpha'(x)|+|\beta'(x)|}{x}+\frac{u^2|\alpha(x)|+|\beta(x)|}{x^2}+\frac{1}{x^4}\right)|T_{2k}(x)|dx.
\end{multline}
Note that various constants are omitted since we are using the $\ll $ sign.
Substituting \eqref{fMellinY4} into \eqref{fMellinY3}, we obtain
\begin{multline}\label{fMellinY5}
\hat{f}_{2k}(0,0;ir)\ll\frac{1}{(1+|r|)^3}\\\times
\int_0^1(1-x)^{2}\left(\frac{u^2|\alpha'(x)|+|\beta'(x)|}{x}+\frac{u^2|\alpha(x)|+|\beta(x)|}{x^2}+\frac{1}{x^4}\right)|T_{2k}(x)|dx\\+
\frac{1}{(1+|r|)^2}\int_0^1(1-x)
\left(\frac{u^2|\alpha(x)|+|\beta(x)|}{x}+\frac{1}{x^3}\right)|T_{2k}(x)|dx.
\end{multline}
Consider the second integral in \eqref{fMellinY5}. Using \eqref{Y def}, \eqref{def fandg} and making the change of variable $x=\cosh^{-2}{\sqrt{\xi}/2}$, we show that
\begin{multline}\label{fMellinY6}
\frac{1}{(1+|r|)^2}\int_0^1(1-x)
\left(\frac{u^2|\alpha(x)|+|\beta(x)|}{x}+\frac{1}{x^3}\right)|T_{2k}(x)|dx\\\ll
\frac{1}{(1+|r|)^2}\int_0^1(1-x)^{1/2}
\left(\frac{u^2}{x^3}+\frac{1}{x^3(1-x)}\right)|\Psi_{2k}(x)|dx\\ \ll
\frac{1}{(1+|r|)^2}
\int_0^{\infty}\Big|\Psi_{2k}\left(\frac{1}{\cosh^2{\sqrt{\xi}/2}} \right)\Big|
\left(u^2\cosh^{2}\frac{\sqrt{\xi}}{2}\sinh^{2}\frac{\sqrt{\xi}}{2}
+\cosh^{4}\frac{\sqrt{\xi}}{2}\right)\frac{d\xi}{\xi^{1/2}}.
\end{multline}
Applying \eqref{psik LG2} and estimating the $K$-Bessel function using \cite[(10.25.3), (10.30.3)]{HMF}, we obtain
\begin{multline}\label{fMellinY7}
\frac{1}{(1+|r|)^2}\int_0^1(1-x)
\left(\frac{u^2|\alpha(x)|+|\beta(x)|}{x}+\frac{1}{x^3}\right)|T_{2k}(x)|dx\\\ll
\frac{1}{(1+|r|)^2}
\int_0^{\infty}\frac{|K_0(u\sqrt{\xi})|}{\sinh^{1/2}{\sqrt{\xi}}}
\left(u^2\cosh^{2}\frac{\sqrt{\xi}}{2}\sinh^{2}\frac{\sqrt{\xi}}{2}
+\cosh^{4}\frac{\sqrt{\xi}}{2}\right)\frac{d\xi}{\xi^{1/4}}\\
\ll\frac{k^{\epsilon}}{k(1+|r|)^2}.
\end{multline}
Consider the first integral in \eqref{fMellinY5}.  Using  \eqref{def fandg}, we have for $0<x<1$
\begin{equation}\label{fandg derivative}
\alpha'(x)\ll\frac{1}{x^3(1-x)}+\frac{1}{x^2(1-x)^2}, \quad
\beta'(x)\ll-\frac{1}{x^3(1-x)^2}+\frac{1}{x^2(1-x)^3}.
\end{equation}
Using  \eqref{def fandg}, \eqref{fandg derivative}, \eqref{Y def} and making the change of variable $x=\cosh^{-2}{\sqrt{\xi}/2}$, we obtain
\begin{multline}\label{fMellinY8}
\frac{1}{(1+|r|)^3}\int_0^1(1-x)^{2}\left(\frac{u^2|\alpha'(x)|+|\beta'(x)|}{x}+\frac{u^2|\alpha(x)|+|\beta(x)|}{x^2}+\frac{1}{x^4}\right)|T_{2k}(x)|dx\\\ll
\frac{1}{(1+|r|)^3}\int_0^1(1-x)^{1/2}
\left(\frac{u^2(1-x)}{x^4}+\frac{u^2}{x^3}+\frac{1}{x^4}+\frac{1}{x^3(1-x)}\right)\left|\Psi_{2k}(x)\right|dx\\\ll 
\frac{1}{(1+|r|)^3}
\int_0^{\infty}\left|\Psi_{2k}\left(\frac{1}{\cosh^2{\sqrt{\xi}/2}} \right)\right|
\Biggl(u^2\cosh^{2}\frac{\sqrt{\xi}}{2}\sinh^{4}\frac{\sqrt{\xi}}{2}+\\
+u^2\cosh^{2}\frac{\sqrt{\xi}}{2}\sinh^{2}\frac{\sqrt{\xi}}{2}
+\cosh^{4}\frac{\sqrt{\xi}}{2}\sinh^{2}\frac{\sqrt{\xi}}{2}
+\cosh^{4}\frac{\sqrt{\xi}}{2}\Biggr)\frac{d\xi}{\xi^{1/2}}.
\end{multline}
Applying \eqref{psik LG2} and the standard bounds on the $K$-Bessel function \cite[(10.25.3), (10.30.3)]{HMF}, we have
\begin{multline}\label{fMellinY9}
\frac{1}{(1+|r|)^3}\int_0^1(1-x)^{2}\left(\frac{u^2|\alpha'(x)|+|\beta'(x)|}{x}+\frac{u^2|\alpha(x)|+|\beta(x)|}{x^2}+\frac{1}{x^4}\right)\\ \times |T_{2k}(x)|dx
\ll\frac{k^{\epsilon}}{k(1+|r|)^2}.
\end{multline}
Substituting \eqref{fMellinY9} and \eqref{fMellinY7} into \eqref{fMellinY5}, we complete the proof of \eqref{fMellin ir}.
\end{proof}

\subsection{Mellin transform of $g_{2k}$}
\begin{lem}\label{lem:gMelTrans}
Assume that $-1/2-\Re{v}<\Re{s}<1/4-\Re{u}/2$ and $0\le\Re{u}<4k-1/2$. Then the Mellin transform of the function $g_{2k}(u,v;x)$ can be written as follows:
\begin{equation}\label{gMellin1}
\hat{g}_{2k}(u,v;s)=
\frac{2^{1/2-u}\sin(\pi(1/4+u/2))}{\pi^{1/2}}
\int_0^1\frac{(1-x)^{s+v-1/2}}{x^{s+u+3/4}}
\Phi_{2k}\left(u;x\right)dx,
\end{equation}
\begin{multline}\label{gMellin2}
\hat{g}_{2k}(u,v;s)=\Gamma(1/2+s+v)
\frac{1}{2\pi i}\int_{(\Delta)}\frac{\Gamma(2k-1/2+w/2)}{\Gamma(2k+1/2-w/2)}\Gamma(\frac{1}{2}-u-w)\\\times
\sin\left( \pi \frac{1/2+u+w}{2}\right)\frac{\Gamma(w/2-s)}{\Gamma(1/2+v+w/2)}2^wdw,
\end{multline}
where $\max(1-4k,2\Re{s})<\Delta<1/2-\Re{u}$,
\begin{multline}\label{gMellin3}
\hat{g}_{2k}(u,v;s)=
\frac{2^{1/2-u}\sin(\pi(1/4+u/2))}{\pi^{1/2}}\Gamma(1/2+s+v)\\\times
\frac{\Gamma(2k-1/4-u/2)\Gamma(3/4-2k-u/2)\Gamma(1/4-u/2-s)}{\Gamma(1/2)\Gamma(3/4+v-u/2)}\\\times
{}_3F_{2}\left(2k-\frac{1}{4}-\frac{u}{2},\frac{3}{4}-2k-\frac{u}{2},\frac{1}{4}-\frac{u}{2}-s;\frac{1}{2},\frac{3}{4}+v-\frac{u}{2}; 1\right).
\end{multline}
\end{lem}
\begin{proof}
It follows from  \eqref{Mellin def} and \eqref{def g} that
\begin{equation}\label{gMellin0}
\hat{g}_{2k}(u,v;s)=
\int_0^{\infty}\frac{x^{s+v-1/2}}{(1+x)^{1/2+v}}I_{2k}\left(u;\frac{2}{(1+x)^{1/2}}\right)dx.
\end{equation}
Applying \eqref{eq:integralI x<2} to evaluate \eqref{gMellin0}, we obtain \eqref{gMellin1}.
Assuming that $\Re{u>0}$, we substitute \eqref{eq:integralI} to \eqref{gMellin0}. For $\Re{u>0},$ $\Re{w}>2\Re{s},$ $\Re{s}>-1/2-\Re{v}$, the resulting double integral converges absolutely. Changing the order of integration and applying \cite[(5.12.3)]{HMF}, namely
\begin{equation*}
\int_0^{\infty}\frac{x^{s+v-1/2}}{(1+x)^{1/2+w/2+v}}dx=
\Gamma(1/2+s+v)\frac{\Gamma(w/2-s)}{\Gamma(1+v+w/2)},
\end{equation*}
we prove \eqref{gMellin2}. Note that the integral on the right-hand side of \eqref{gMellin2} converges absolutely provided that $\Re{s}>-1/2-\Re{u}-\Re{v}.$ Moving the line of integration in \eqref{gMellin2} to the right and crossing the poles at $w=1/2-u+j$, we obtain \eqref{gMellin3}.

\end{proof}

\begin{lem}
For $\Re{v}>-1/4$ the following equality holds
\begin{equation}\label{gMellin -1/4}
\hat{g}_{2k}(0,v;-1/4)=-\sqrt{2}\sin(\pi v)
\frac{\Gamma^2(1/4+v)}{\pi^{1/2}}
\frac{\Gamma(2k-1/4)\Gamma(2k-v)}{\Gamma(2k+1/4)\Gamma(2k+v)}.
\end{equation}
In particular,
\begin{equation}\label{gMellin -1/4 v=0}
\hat{g}_{2k}(0,0;-1/4)=0.
\end{equation}
\end{lem}
\begin{proof}
According to \eqref{gMellin3} we have
\begin{multline*}
\hat{g}_{2k}(0,v;-1/4)=
\frac{\Gamma(1/4+v)\Gamma(2k-1/4)\Gamma(3/4-2k)}{\Gamma(3/4+v)\pi^{1/2}}
{}_2F_{1}\left(2k-\frac{1}{4},\frac{3}{4}-2k;\frac{3}{4}+v; 1\right).
\end{multline*}
Applying \cite[(15.4.20)]{HMF}, this expression simplifies to
\begin{equation*}
\hat{g}_{2k}(0,v;-1/4)=
\frac{\Gamma^2(1/4+v)}{\pi^{1/2}}
\frac{\Gamma(2k-1/4)\Gamma(3/4-2k)}{\Gamma(1+v-2k)\Gamma(2k+v)}.
\end{equation*}
Finally, using \cite[(5.5.3)]{HMF} we obtain \eqref{gMellin -1/4}.
\end{proof}

\begin{lem}\label{Lem gMellin v 1/2-v}
For $v\to 1/4$ the following asymptotic formulas hold
\begin{equation}\label{gMellin v 1/2-v}
\hat{g}_{2k}(0,v;1/2-v)=\frac{2^{3/2}}{2v-1/2}\frac{\Gamma(2k-1/4)}{\Gamma(2k+1/4)}+O(1),
\end{equation}
\begin{equation}\label{gMellin dif v 1/2-v}
\frac{\partial}{\partial s}\hat{g}_{2k}(0,v;s)\Bigg|_{s=1/2-v}=
\frac{2^{5/2}}{(2v-1/2)^2}\frac{\Gamma(2k-1/4)}{\Gamma(2k+1/4)}+
O(1).
\end{equation}
Furthermore,
\begin{equation}\label{gMellin 1/2}
\hat{g}_{2k}(0,0;1/2)=2^{3/2}\Gamma(-1/2)+O(k^{-1+\epsilon}),
\end{equation}
\begin{multline}\label{gMellin dif 1/2}
\frac{\partial}{\partial s}\hat{g}_{2k}(0,0;s)\Bigg|_{s=1/2}=
-2^{5/2}\Gamma(-1/2)\frac{\Gamma'}{\Gamma}(-1/2)\\+
2^{3/2}\Gamma(-1/2)\left(2\frac{\Gamma'}{\Gamma}(2k)+2\log2-\pi\right)+O(k^{-1+\epsilon}).
\end{multline}
\end{lem}
\begin{proof}
For $u=0$, $0\le\Re{v}\le1/2$  and $0<\Re{s}<1/4$, we move the line of integration in \eqref{gMellin2} to $-2+2\Re{s}<\Delta<2\Re{s}$ crossing the pole at $w=2s$. Hence
\begin{multline}\label{gMellin 0vs}
\hat{g}_{2k}(0,v;s)=
2^{2s+1}\frac{\Gamma(2k-1/2+s)}{\Gamma(2k+1/2-s)}\Gamma(\frac{1}{2}-2s)\sin\left( \pi \frac{1/2+2s}{2}\right)\\+
\Gamma(1/2+s+v)\frac{1}{2\pi i}\int_{(\Delta)}\frac{\Gamma(2k-1/2+w/2)}{\Gamma(2k+1/2-w/2)}\Gamma(\frac{1}{2}-u-w)\\ \times
\sin\left( \pi \frac{1/2+u+w}{2}\right)\frac{\Gamma(w/2-s)}{\Gamma(1/2+v+w/2)}2^wdw,
\end{multline}
where $-2+2\Re{s}<\Delta<\min(2\Re{s},1/2).$  Therefore, \eqref{gMellin 0vs} is now valid for $\Re{s}<5/4.$
Choosing $\Delta=0$, we obtain
\begin{multline}\label{gMellin 0v1/2-v 2}
\hat{g}_{2k}(0,v;1/2-v)=
2^{2-2v}\frac{\Gamma(2k-v)}{\Gamma(2k+v)}\Gamma(2v-\frac{1}{2})\sin(3\pi/4-\pi v)+\\+
\frac{1}{\pi}\int_{-\infty}^{\infty}\frac{\Gamma(2k-1/2+ir)}{\Gamma(2k+1/2-ir)}\Gamma(\frac{1}{2}-2ir)
\frac{\sin(\pi/4+\pi ir)2^{2ir}dr}{(ir-1/2+v)}.
\end{multline}
Estimating the integral above trivially using Stirling's formula, we conclude the proof of \eqref{gMellin 1/2}.
Another direct consequence of the representation \eqref{gMellin 0v1/2-v 2} is \eqref{gMellin v 1/2-v}.
Finally, the formulas \eqref{gMellin dif v 1/2-v} and \eqref{gMellin dif 1/2} can also be derived from \eqref{gMellin 0vs} by taking the derivative with respect to $s$.
\end{proof}


\begin{lem}\label{Lem gMellin big r}
For $|r|>3k$ and any $A>0$ we have
\begin{equation}\label{gMellin big r}
\hat{g}_{2k}(0,0;ir)\ll\frac{1}{|r|^A}.
\end{equation}
\end{lem}
\begin{proof}
Using the representation \eqref{gMellin3}, we obtain
\begin{equation*}
\hat{g}_{2k}(0,0;ir)=
\frac{\Gamma(1/2+ir)}{\pi^{1/2}}
\sum_{j=0}^{\infty}\frac{(-1)^j}{j!}
\frac{\Gamma(2k-1/4+j)\Gamma(3/4-2k+j)\Gamma(1/4-ir+j)}{\Gamma(1/2+j)\Gamma(3/4+j)}.
\end{equation*}
It follows from \cite[(5.4.4), (5.6.6)]{HMF} that 
\begin{equation*}|\Gamma(1/2+ir)|\ll\exp(-\pi|r|/2)\text{ and }|\Gamma(1/4-ir+j)|\ll\Gamma(1/4+j).
\end{equation*}
 Consequently,
\begin{equation*}
\hat{g}_{2k}(0,0;ir)\ll
\exp(-\pi|r|/2)
\sum_{j=0}^{\infty}
\frac{|\Gamma(3/4-2k+j)|\Gamma(2k-1/4+j)\Gamma(1/4+j)}{\Gamma(1+j)\Gamma(1/2+j)\Gamma(3/4+j)}.
\end{equation*}
According to \cite[(5.5.3)]{HMF} we have $|\Gamma(3/4-2k+j)|=|\Gamma(2k+1/4-j)|^{-1}.$ Furthermore,
\begin{equation*}
\frac{\Gamma(1/4+j)}{\Gamma(1+j)\Gamma(1/2+j)\Gamma(3/4+j)}\ll\frac{1}{\Gamma^2(1+j)}.
\end{equation*}
As a result,
\begin{multline}\label{gMellin big r1}
\hat{g}_{2k}(0,0;ir)\ll\exp(-\pi|r|/2)
\sum_{j=0}^{2k-1}
\frac{\Gamma(2k-1/4+j)}{\Gamma(2k+1/4-j)\Gamma^2(1+j)}+\\
\exp(-\pi|r|/2)\sum_{j=2k}^{\infty}
\frac{\Gamma(2k-1/4+j)\Gamma(j+3/4-2k)}{\Gamma^2(1+j)}.
\end{multline}
Using Stirling's formula \cite[(5.11.3)]{HMF} we obtain that for $0<j<2k-1$
\begin{multline*}
\frac{\Gamma(2k-1/4+j)}{\Gamma(2k+1/4-j)\Gamma^2(1+j)}\asymp
\frac{(2k+j)^{-1/4}}{(2k-j)^{1/4}j^2}\frac{\Gamma(2k+j)}{\Gamma(2k-j)\Gamma^2(j)}\asymp\\
\frac{(2k+j)^{-3/4}}{(2k-j)^{-1/4}j}
\frac{(2k+j)^{2k+j}}{(2k-j)^{2k-j}j^{2j}}=
\frac{(2k-j)^{1/4}}{(2k+j)^{3/4}j}\exp(s_1(k,j)),
\end{multline*}
where
\begin{equation*}
s_1(k,j):=(2k+j)\log(2k+j)-(2k-j)\log(2k-j)-2j\log j.
\end{equation*}
The function $s_1(k,j)$ attains its maximum at the point $j=k\sqrt{2}$, and therefore,
\begin{multline}\label{gMellin big r2}
\exp(-\pi|r|/2)\sum_{j=0}^{2k-1}
\frac{\Gamma(2k-1/4+j)}{\Gamma(2k+1/4-j)\Gamma^2(1+j)}\ll\\
\exp(-\pi|r|/2)\sum_{j=1}^{2k-1}\frac{(2k-j)^{1/4}}{(2k+j)^{3/4}j}\exp(s_1(k,j))\ll\\
\exp(-\pi|r|/2+s_1(k,k\sqrt{2}))k^{1/2}\ll\frac{1}{|r|^A}
\end{multline}
for $|r|>3k.$ In the same way  we show that for $j>2k$
\begin{equation*}
\frac{\Gamma(2k-1/4+j)\Gamma(j+3/4-2k)}{\Gamma^2(1+j)}\asymp
\frac{(j-2k)^{1/4}}{(2k+j)^{3/4}j}\exp(s_2(k,j)),
\end{equation*}
where
\begin{equation*}
s_2(k,j):=(2k+j)\log(2k+j)-(j-2k)\log(j-2k)-2j\log j.
\end{equation*}
The function $s_2(k,j)$ is decreasing. For $j=2k$ it follows from \cite[(5.5.5)]{HMF} that
\begin{equation*}
\frac{\Gamma(2k-1/4+j)\Gamma(j+3/4-2k)}{\Gamma^2(1+j)}\asymp
\frac{\Gamma(4k-1/4)}{\Gamma^2(2k+1)}\ll\frac{2^{4k}}{k^{7/4}}.
\end{equation*}
Finally, we obtain
\begin{multline}\label{gMellin big r3}
\exp(-\pi|r|/2)\sum_{j=2k}^{\infty}
\frac{\Gamma(2k-1/4+j)\Gamma(j+3/4-2k)}{\Gamma^2(1+j)}\ll\\
\exp(-\pi|r|/2)\sum_{j=2k}^{\infty}\frac{2^{4k}}{j^{3/2}}\ll\frac{1}{|r|^A}
\end{multline}
for $|r|>3k.$ Substituting \eqref{gMellin big r2} and \eqref{gMellin big r3} into \eqref{gMellin big r1}, we prove the lemma.
\end{proof}

To investigate the behavior of the function  $\hat{g}_{2k}(0,0;ir)$ for $|r|\le3k$  we use the formula \eqref{gMellin1}.

\begin{lem}\label{lem:LGapprox}
Let $\kc:=4k-1.$ We have
\begin{multline}\label{gMellin BessInt}
\hat{g}_{2k}(0,0;ir)=-2^{3/2}\pi^{1/2}
\int_0^{\pi/2}\frac{(\tan x)^{2ir}}{(\sin(2x))^{1/2}}
Y_0(\kc x)x^{1/2}dx-\\
-2^{3/2}\pi^{1/2}
\int_0^{\pi/2}\frac{(\tan x)^{2ir}}{(\sin(2x))^{1/2}}
J_0(\kc x)x^{1/2}dx+O(k^{-3/2}),
\end{multline}
\begin{equation}\label{gMellin triv.est.}
\hat{g}_{2k}(0,0;ir)\ll\frac{1}{k^{1/2}}.
\end{equation}
\end{lem}
\begin{proof}
It follows from \eqref{gMellin1} that
\begin{equation}\label{gMellin1 at 00}
\hat{g}_{2k}(0,0;ir)=
\frac{1}{\pi^{1/2}}
\int_0^{\pi^2/4}\frac{(\tan\sqrt{\xi})^{2ir}}{(\cos\sqrt{\xi})^{1/2}}
\Phi_{2k}\left(\cos^2\sqrt{\xi}\right)\frac{d\xi}{\sqrt{\xi}}.
\end{equation}
In order to prove \eqref{gMellin BessInt}, we apply the following approximation of the function $\Phi_{2k}$
(see  \cite[Theorems~6.5 and 6.10, Corollary~6.9]{BF1} for details):
\begin{multline}\label{eq:Phi2k}
\Phi_{2k}(\cos^2\sqrt{\xi}) =
\frac{-\pi}{\xi^{1/4}(\sin{\sqrt{\xi}})^{1/2}}
\Biggl[\sqrt{\xi}Y_0((4k-1)\sqrt{\xi}) + \\
\sqrt{\xi}J_0((4k-1)\sqrt{\xi})
+ O\left(\frac{1}{k}\left|\sqrt{\xi}Y_0((4k-1)\sqrt{\xi})\right|\right)\Biggr].
\end{multline}
Substituting \eqref{eq:Phi2k} into \eqref{gMellin1 at 00}  and estimating the error term using standard estimates on the $Y$-Bessel function \cite[(10.7.1), (10.7.8)]{HMF}, we obtain \eqref{gMellin BessInt}. Applying \cite[(10.7.1), (10.7.8)]{HMF} to estimate the integrals in \eqref{gMellin BessInt} we prove \eqref{gMellin triv.est.}.
\end{proof}

The estimate \eqref{gMellin triv.est.} is sufficiently good for our purposes only if $r\ll \kc^{1/2-\delta}$. For $r\gg\kc^{1/2-\delta}$, it is required to analyse \eqref{gMellin BessInt} more carefully. We consider further only the first integral  in \eqref{gMellin BessInt}, as the second integral can be treated similarly.
The idea is to replace the $Y$-Bessel function  in \eqref{gMellin BessInt} by its asymptotic formula  \cite[(10.7.8)]{HMF}. To this end, we first make the following partition of unity:
\begin{equation}\label{partition of unity}
\alpha_1(x)+\beta(x)+\alpha_2(x)=1 \quad\hbox{for}\quad0\le x\le\frac{\pi}{2},
\end{equation}
where $\alpha_{1,2}(x)$, $\beta(x)$ are smooth infinitely differentiable functions such that for some small $\varepsilon>0$ (to be chosen later), we have
\begin{equation}\label{partition of unity1}
\alpha_1(x)=1 \text{ for } 0\le x\le \varepsilon,\quad
\alpha_1(x)=0 \text{ for } x\ge 2\varepsilon,
\end{equation}
\begin{equation}\label{partition of unity2}
\alpha_2(x)=1 \text{ for } \frac{\pi}{2}-\varepsilon\le x\le\frac{\pi}{2},\quad
\alpha_2(x)=0 \text{ for } 0\le x\le \frac{\pi}{2}-2\varepsilon,
\end{equation}
\begin{equation}\label{partition of unity3}
\beta(x)=1 \text{ for }2\varepsilon\le x\le \frac{\pi}{2}-2\varepsilon,\quad 
\beta(x)=0 \text{ for } 0\le x\le \varepsilon, \frac{\pi}{2}-\varepsilon\le x\le \frac{\pi}{2},
\end{equation}
and $\alpha_{1,2}^{(j)}(x)\ll\varepsilon^{-j}$, $\beta^{(j)}(x)\ll\varepsilon^{-j}.$
\begin{lem}
For $|r|>1$  the following holds
\begin{multline}\label{gMellin bettaBessInt}
\hat{g}_{2k}(0,0;ir)\ll
\left|\int_0^{\pi/2}\frac{\beta(x)(\tan x)^{2ir}}{(\sin(2x))^{1/2}}
Y_0(\kc x)x^{1/2}dx\right|+
\\
\left|
\int_0^{\pi/2}\frac{\beta(x)(\tan x)^{2ir}}{(\sin(2x))^{1/2}}
J_0(\kc x)x^{1/2}dx\right|
+O\left(\frac{k^{-1+\epsilon}+k^{1/2}\varepsilon^{3/2}}{r}\right).
\end{multline}
\end{lem}
\begin{proof}
As the first step, we use the partition of unity \eqref{partition of unity} to rewrite the integrals in \eqref{gMellin BessInt}. Then to prove the lemma, it is required to estimate the contribution of integrals with $\alpha_{1,2}(x).$ All these integrals can be analysed similarly. Therefore, we consider only
\begin{equation}\label{Y alpha 1}
I_1:=\int_0^{\pi/2}\frac{\alpha_1(x)Y_0(\kc x)x^{1/2}}{(\sin(2x))^{1/2}}
(\tan x)^{2ir}dx.
\end{equation}
Integrating by parts we obtain
\begin{equation}\label{Y alpha 2}
I_1\ll\frac{1}{r}\int_0^{\pi/2}\frac{\partial}{\partial x}\left[
\alpha_1(x)Y_0(\kc x)x^{1/2}(\sin(2x))^{1/2}\right](\tan x)^{2ir}dx.
\end{equation}
Evaluating the derivative and estimating the integral trivially with the use of \cite[(10.7.1), (10.7.8)]{HMF}, we complete the proof of \eqref{gMellin bettaBessInt}.
\end{proof}
For simplicity, let us assume further that $r>0$. The case $r<0$ can be treated in the same way.

\begin{lem}\label{int:saddle}
Let $\delta$ be some fixed constant such that $0<\delta<1/4$. Then  for $r>1$ and $\varepsilon=k^{-1/2-2\delta}$ we have
\begin{equation}\label{gMellin bettaInt}
\hat{g}_{2k}(0,0;ir)\ll\frac{1}{k^{1/2}}
\left|\int_0^{\pi/2}\frac{\beta(x)\exp(i\kc h(x))}{(\sin(2x))^{1/2}}
dx\right|+\frac{k^{-1/4-3\delta}+k^{-1/2}}{r}+\frac{1}{k^{5/4-\delta}},
\end{equation}
where
\begin{equation}\label{h def}
h(x)=-x+\frac{2r}{\kc}\log(\tan x).
\end{equation}
\end{lem}
\begin{proof}
We substitute the asymptotic formulas for Bessel functions \cite[(8.451.1), (8.451.2)]{GR} into \eqref{gMellin bettaBessInt} and estimate the error terms by its absolute value, obtaining
\begin{equation*}
\hat{g}_{2k}(0,0;ir)\ll\frac{1}{k^{1/2}}
\sum_{\pm}\left|\int_0^{\pi/2}\frac{\beta(x)\exp(ih_{\pm}(x))}{(\sin(2x))^{1/2}}
dx\right|+\frac{k^{-1/4-3\delta}}{r}+\frac{1}{k^{5/4-\delta}},
\end{equation*}
where
\begin{equation}\label{h_pm def}
h_{\pm}(x)=\pm\kc x+2r\log(\tan x).
\end{equation}
The integral with $h_{+}(x)$ can be estimated using \cite[Lemma 4.3]{T}. More precisely,
\begin{equation*}
\frac{1}{k^{1/2}}
\left|\int_0^{\pi/2}\frac{\beta(x)\exp(ih_{+}(x))}{(\sin(2x))^{1/2}}
dx\right|\ll\frac{1}{k^{1/2}}\max_{0<x<\pi/2}\frac{\beta(x)(\sin(2x))^{1/2}}{\kc\sin(2x)+4r}\ll\frac{1}{rk^{1/2}}.
\end{equation*}
\end{proof}


The classical approach to estimate the integral on the right-hand side of \eqref{gMellin bettaInt} is the saddle point method (also called the method of steepest descent). Another possibility is the stationary phase method, which is in some sense (see discussion in \cite[pp. 276--279]{BlHan}) an analogue of the saddle point method for Fourier-type integrals. The first step in all these methods is to determine the so-called saddle points of the function $h(x)$ defined as zeros of $h'(x)=0.$ Using \eqref{h def} we find that
\begin{equation}\label{h derivative}
h'(x)=-1+\frac{4r}{\kc}\frac{1}{\sin(2x)}.
\end{equation}
It is convenient to introduce two new parametrs $\vartheta$ and $\mu$ such that:
\begin{equation}\label{varteta def}
\sin(2\vartheta)=\frac{4r}{\kc},\quad0<\vartheta<\frac{\pi}{4}\quad\hbox{if}\quad 4r\le\kc,
\end{equation}
\begin{equation}\label{mu def}
\cosh(2\mu)=\frac{4r}{\kc},\quad \mu>0\quad\hbox{if}\quad 4r>\kc.
\end{equation}
Then the saddle points of the function $h(x)$ are
\begin{equation}\label{sadpoints def1}
x_1=\vartheta, \quad x_2=\frac{\pi}{2}-\vartheta\quad\hbox{if}\quad 4r\le\kc,
\end{equation}
\begin{equation}\label{sadpoints def2}
x_3=\frac{\pi}{4}-i\mu, \quad x_4=\frac{\pi}{4}+i\mu\quad\hbox{if}\quad 4r>\kc.
\end{equation}
We consider only the case $4r\le\kc$ since the second case can be analysed similarly.
Note that the condition $r\gg\kc^{1/2-\delta}$ implies that
\begin{equation}\label{teta conditions}
\vartheta\gg\kc^{-1/2-\delta}>k^{-1/2-2\delta}=\varepsilon.
\end{equation}
Thus both saddle points belong to the interval of integration.

An important observation is that as  $4r\rightarrow\kc$ {\it the saddle points coalesce}. It is known that in this case the behaviour of the integral changes. Therefore, it is required to analyse three different ranges: 
\begin{itemize} \item $r$ is small, \item $r$ is near $\kc/4$, \item $r$ is large.
\end{itemize} 

The case of coalescing saddle points is usually described in books, see  \cite[Section 9.2]{BlHan} and \cite[Section 7.4]{Wong}. It is well known that the standard saddle point method does not work in this situation and a more refined analysis is required. Therefore, we mainly follow \cite[Section 9.2]{BlHan}. This approach was originally developed by Chester, Friedman and Ursell \cite{CFU}, with some additional ideas due to Bleistein \cite{Bl}.


The main idea of the method is to change the variable of integration such that the integral can be written in terms of the Airy function, which has for real $x$ the following representation:
\begin{equation}\label{Airy def}
Ai(ax^{2/3})=\frac{x^{1/3}}{2\pi}\int_{-\infty}^{\infty}\exp\left(ix\left(\frac{y^3}{3}+ay\right)\right)dy.
\end{equation}
For simplicity, let us denote
\begin{equation}\label{g def}
g(x):=\frac{\beta(x)}{(\sin(2x))^{1/2}}.
\end{equation}
Our goal is to estimate the integral (see \eqref{gMellin bettaInt})
\begin{equation}\label{I Int}
I=\int_0^{\pi/2}g(x)\exp(i\kc h(x))dx, \quad h(x)=-x+\frac{\sin\vartheta}{2}\log(\tan x).
\end{equation}
To this end, following \cite[(9.2.6)]{BlHan} we define a new variable $t$ such that:
\begin{equation}\label{t variable def}
h(x)=\frac{t^3}{3}-\gamma^2 t+\rho,
\end{equation}
where the constants $\gamma$ and $\rho$ are chosen such that the point $x=x_1$ corresponds to $t=-\gamma$ and  the point $x=x_2$ corresponds to $t=\gamma.$  These conditions yield
\begin{equation*}
h(x_1)=\frac{2\gamma^3}{3}+\rho,\quad h(x_2)=-\frac{2\gamma^3}{3}+\rho,
\end{equation*}
and therefore,
\begin{equation}\label{gamma rho expression1}
\frac{4\gamma^3}{3}=h(x_1)-h(x_2),\quad 2\rho=h(x_1)+h(x_2).
\end{equation}
Evaluating $h(x_{1,2})$ we find that $\rho=-\pi/4$ and
\begin{equation}\label{gamma rho expression2}
\frac{4\gamma^3}{3}=\frac{\pi}{2}-2\vartheta+\sin(2\vartheta)\log(\tan\vartheta).
\end{equation}
Note that for $0<\vartheta<\pi/4$ the right-hand side of \eqref{gamma rho expression2} is positive, and that  for $\vartheta=\pi/4$ we obtain $\gamma=0$.  Changing
the variable $x$ in the integral \eqref{I Int} by $t$ defined by \eqref{t variable def}, we obtain an analogue of \cite[(9.2.18), (9.2.19)]{BlHan}, namely
\begin{equation}\label{I Int2}
I=\exp(i\kc \rho)\int_{-\infty}^{\infty}G_0(t,\vartheta)\exp\left(i\kc (t^3/3-\gamma^2 t)\right)dt,
\end{equation}
where
\begin{equation}\label{G0 def}
G_0(t,\vartheta)=g(x(t))\frac{dx}{dt}.
\end{equation}
As in \cite[(9.2.20)]{BlHan} we write
\begin{equation}\label{H0 def}
G_0(t,\vartheta)=a_0+a_1t+(t^2-\gamma^2)H_0(t,\vartheta),
\end{equation}
where (see \cite[(9.2.21),(9.2.22)]{BlHan})
\begin{equation}\label{a0a1 def}
a_0=\frac{G_0(\gamma,\vartheta)+G_0(-\gamma,\vartheta)}{2},\quad
a_1=\frac{G_0(\gamma,\vartheta)-G_0(-\gamma,\vartheta)}{2\gamma}.
\end{equation}
Note that  $a_0$ and $a_1$ are chosen such that the function $H_0(t,\vartheta)$ is regular at points $t=\pm\gamma.$

To evaluate $a_{0,1}$, as well as to analyse the properties of $H_0(t,\vartheta)$, we need some preliminary results.
\begin{lem}\label{lem derivatives}
For $\vartheta<\pi/4$ we have
\begin{equation}\label{dx/dt at gamma}
\frac{dx}{dt}\Biggr|_{t=\pm\gamma}=\sqrt{\gamma\tan(2\vartheta)},
\end{equation}
\begin{equation}\label{d2x/dt2 at gamma}
\frac{d^2x}{dt^2}\Biggr|_{t=\pm\gamma}=\mp\frac{1}{3}\left(
4\gamma+2\gamma\tan^2(2\vartheta)-\sqrt{\frac{\tan(2\vartheta)}{\gamma}}
\right).
\end{equation}
For $\vartheta=\pi/4$ we have
\begin{equation}\label{dx/dt at 0}
\frac{dx}{dt}\Biggr|_{t=0}=2^{-1/3},\quad
\frac{d^2x}{dt^2}\Biggr|_{t=0}=0,\quad
\frac{d^3x}{dt^3}\Biggr|_{t=0}=-1.
\end{equation}
\end{lem}
\begin{proof}
First, consider the case $\vartheta<\pi/4$, $t=-\gamma$. We can write
\begin{equation}\label{taylor series 1}
x-\vartheta=\sum_{n=0}^{\infty}b_n(t+\gamma)^n,\quad
h'(x)=\sum_{n=0}^{\infty}c_n(t+\gamma)^n.
\end{equation}
Let us compute $b_i$, $c_i$ for $i=0,1,2.$ Note that $b_0=0$ since the point $x=\vartheta$ corresponds to $t=-\gamma$. We have
\begin{equation}\label{taylor series 2}
h'(x)=-1+\frac{\sin(2\vartheta)}{\sin(2x)}=-\frac{2(x-\vartheta)}{\tan(2\vartheta)}+
(x-\vartheta)^2\left(\frac{4}{\tan^2(2\vartheta)}+2\right)+O((x-\vartheta)^3).
\end{equation}
Substituting the expansion for $(x-\vartheta)$ from \eqref{taylor series 1} into \eqref{taylor series 2}, we show that
\begin{equation}\label{c coeff}
c_0=0,\quad
c_1=-\frac{2b_1}{\tan(2\vartheta)},\quad
c_2=-\frac{2b_2}{\tan(2\vartheta)}+b_1^2\left(\frac{4}{\tan^2(2\vartheta)}+2\right).
\end{equation}
It follows from \eqref{t variable def} that
\begin{equation}\label{t variable def conseq}
h'(x)\frac{dx}{dt}=t^2-\gamma^2=-2\gamma(t+\gamma)+(t+\gamma)^2.
\end{equation}
Substituting \eqref{taylor series 1} into \eqref{t variable def conseq} yields
\begin{equation}\label{b,c coeff relation}
c_1b_1=-2\gamma,\quad
c_2b_1+2c_1b_2=1.
\end{equation}
Using \eqref{b,c coeff relation} and \eqref{c coeff}, we obtain
\begin{equation*}
b_1=\sqrt{\gamma\tan(2\vartheta)},\quad
b_2=\frac{1}{6}\left(
4\gamma+2\gamma\tan^2(2\vartheta)-\sqrt{\frac{\tan(2\vartheta)}{\gamma}}
\right).
\end{equation*}
This proves \eqref{dx/dt at gamma} and \eqref{d2x/dt2 at gamma} for $t=-\gamma$. The case $t=\gamma$ is similar.

Second, consider  $\vartheta=\pi/4$. In that case $\gamma=0.$ We can write
\begin{equation}\label{taylor series 3}
x-\frac{\pi}{4}=\sum_{n=0}^{\infty}d_nt^n,\quad
h'(x)=\sum_{n=0}^{\infty}e_nt^n.
\end{equation}
We proceed to compute $d_i$, $e_i$ for $i=0,1,2,3.$ Note that $d_0=0$. Furthermore, we have
\begin{equation}\label{taylor series 4}
h'(x)=-1+\frac{1}{\sin(2x)}=2\left(x-\frac{\pi}{4}\right)^2+\frac{10}{3}\left(x-\frac{\pi}{4}\right)^4
+O((x-\pi/4)^6).
\end{equation}
Substituting the expansion for $(x-\pi/4)$ from \eqref{taylor series 3} into \eqref{taylor series 4}, we show that
\begin{equation}\label{e coeff}
e_0=e_1=0,\quad
e_2=2d_1^2,\quad
e_3=4d_1d_2,\quad
e_4=2d^2_2+4d_1d_3+\frac{10}{3}d_1^4.
\end{equation}
Substituting \eqref{taylor series 3} into \eqref{t variable def conseq} gives
\begin{equation}\label{d,e coeff relation}
e_2d_1=1,\quad
2e_2d_2+e_3d_1=0,\quad
3e_2d_3+2e_3d_2+e_4d_1=0.
\end{equation}
Using \eqref{d,e coeff relation} and \eqref{e coeff} we finally show that
\begin{equation*}
d_1=2^{-1/3},\quad
d_2=0,\quad
d_3=-\frac{1}{6}.
\end{equation*}
This completes the proof of \eqref{dx/dt at 0}.
\end{proof}

\begin{lem}\label{lem a0a1}
For $\vartheta<\pi/4$ we have
\begin{equation}\label{a0a1 value}
a_0=\sqrt{\frac{\gamma}{\cos(2\vartheta)}},\quad
a_1=0,
\end{equation}
and for $\vartheta=\pi/4$ we have
\begin{equation}\label{a0a1 pi/4value}
a_0=2^{-1/3},\quad a_1=0.
\end{equation}
\end{lem}
\begin{proof}
Consider the case $\vartheta<\pi/4$. It follows from \eqref{a0a1 def}, \eqref{G0 def}, \eqref{dx/dt at gamma} and  \eqref{g def} that
\begin{equation*}
a_0=\frac{\beta(\vartheta)+\beta(\pi/2-\vartheta)}{2}\sqrt{\frac{\gamma}{\cos(2\vartheta)}},\quad
a_1=\frac{\beta(\pi/2-\vartheta)-\beta(\vartheta)}{2\gamma}\sqrt{\frac{\gamma}{\cos(2\vartheta)}}.
\end{equation*}
As a consequence of \eqref{teta conditions} we obtain \eqref{a0a1 value}.

Consider the case $\vartheta=\pi/4$. It follows from \eqref{a0a1 def}, \eqref{G0 def}, \eqref{dx/dt at 0} and  \eqref{g def} that $a_0=2^{-1/3}$ and
\begin{equation*}
a_1=\frac{d}{dt}G_0\left(t,\frac{\pi}{4}\right)\Bigr|_{t=0}=
g'(\pi/4)\left(\frac{dx}{dt}\Bigr|_{t=0}\right)^2+g(\pi/4)\frac{d^2x}{dt^2}\Bigr|_{t=0}=
\frac{d^2x}{dt^2}\Bigr|_{t=0}=0.
\end{equation*}
This proves \eqref{a0a1 pi/4value}.
\end{proof}

Substituting \eqref{H0 def} into \eqref{I Int2} and using Lemma \ref{lem a0a1}, we obtain the following representation for our integral 
\begin{multline*}
I=\exp(i\kc \rho)a_0\int_{-\infty}^{\infty}\exp\left(i\kc (t^3/3-\gamma^2 t)\right)dt+\\
\frac{\exp(i\kc \rho)}{i\kc}\int_{-\infty}^{\infty}H_0(t,\vartheta)d\exp\left(i\kc (t^3/3-\gamma^2 t)\right).
\end{multline*}
Using \eqref{Airy def} and integrating by parts yields
\begin{equation*}
I=\exp(i\kc \rho)\frac{2a_0\pi}{\kc^{1/3}}Ai(-\gamma^2\kc^{2/3})
-\frac{\exp(i\kc \rho)}{i\kc}\int_{-\infty}^{\infty}\frac{d}{dt}\left(H_0(t,\vartheta)\right)\exp\left(i\kc (t^3/3-\gamma^2 t)\right)dt.
\end{equation*}
Since the function $H_0(t,\vartheta)$ has a finite number of intervals of monotonicity, we can  estimate the integral in the formula above by its absolute value, getting
\begin{equation}\label{I Int4}
I=\exp(i\kc \rho)\frac{2a_0\pi}{\kc^{1/3}}Ai(-\gamma^2\kc^{2/3})+
O\left(\frac{1}{\kc}\max_{t}|H_0(t,\vartheta)|\right).
\end{equation}
The final step is to estimate the function $H_0(t,\vartheta)$.

\begin{lem}
For $\vartheta\le\pi/4$ satisfying \eqref{teta conditions}  we have
\begin{equation}\label{H0 estimate}
\max_{t}|H_0(t,\vartheta)|\ll\frac{1}{\sqrt{\vartheta}}.
\end{equation}
\end{lem}
\begin{proof}
It is required to consider separately three different cases: $\vartheta\rightarrow\frac{\pi}{4},$ $\vartheta\rightarrow0$ and the third case when $\vartheta$ is some fixed number.

First, let us assume that $\vartheta$ is some fixed number. Then $\gamma$ (see \eqref{gamma rho expression2}) is also some fixed number. We start by estimating $H_0(t,\vartheta)$ near the points $t=\pm\gamma.$ According to \cite[(9.2.24)]{BlHan}, we have
\begin{equation}\label{H0 at pm gamma}
\lim_{t\rightarrow\pm\gamma}H_0(t,\vartheta)=\pm\frac{1}{2\gamma}\frac{d}{dt}G_0(t,\vartheta)\Bigr|_{t=\pm\gamma}.
\end{equation}
Using \eqref{g def}, \eqref{G0 def} and Lemma \ref{lem derivatives}, we prove that $H_0(t,\vartheta)$ is bounded near the points $t=\pm\gamma.$ Other critical points of $H_0(t,\vartheta)$ are $t\rightarrow\pm\infty.$ In this case we use \eqref{g def}, \eqref{H0 def} and \eqref{t variable def conseq}, getting
\begin{equation}\label{H0 expression}
H_0(t,\vartheta)=\frac{g(x(t))}{h'(x(t))}-\frac{a_0}{t^2-\gamma^2}=
\frac{\beta(x(t))\sqrt{\sin(2x(t))}}{\sin(2\vartheta)-\sin(2x(t))}-\frac{a_0}{t^2-\gamma^2}.
\end{equation}
Consequently, for all $t$ that do not belong to a neighbourhood of $\pm\gamma$, the function $H_0(t,\vartheta)$ is trivially bounded by a constant.

Second, consider the case $\vartheta\rightarrow\frac{\pi}{4}.$ It is enough to prove \eqref{H0 estimate} for $\vartheta=\pi/4.$ In this case we have $\gamma=0$ and (see \eqref{H0 at pm gamma})
\begin{equation}\label{H0 at 0}
\lim_{t\rightarrow0}H_0(t,\pi/4)=\frac{1}{2}\frac{d^2}{dt^2}G_0(t,\pi/4)\Bigr|_{t=0}.
\end{equation}
It follows from \eqref{G0 def} that
\begin{equation}\label{G0 second der}
\frac{d^2}{dt^2}G_0(t,\pi/4)=
g''(x(t))\left(\frac{dx}{dt}\right)^3+3g'(x(t))\frac{dx}{dt}\frac{d^2x}{dt^2}+
g(x(t))\frac{d^3x}{dt^3}.
\end{equation}
Using the fact that all derivatives of $g(x)$ at $x=\pi/4$  are finite and applying Lemma \ref{lem derivatives}, we conclude that the limit in \eqref{H0 at 0} is finite.
For $t$ outside of a neighbourhood of $0$, the function $H_0(t,\pi/4)$ is trivially bounded by a constant using \eqref{H0 expression}.

Third,  consider the case $\vartheta\rightarrow0$.  Let us estimate the right-hand side of \eqref{H0 at pm gamma}. We remark that $\gamma$ is a constant for small $\vartheta$ . Therefore, it is only required to estimate the derivative of $G_0(t,\vartheta).$ It follows from \eqref{G0 def} that
\begin{equation*}
\frac{d}{dt}G_0(t,\vartheta)=g'(x(t))\left(\frac{dx}{dt}\right)^2+g(x(t))\frac{d^2x}{dt^2}.
\end{equation*}
Using \eqref{g def} and  Lemma \ref{lem derivatives} we obtain the estimate
\begin{equation}\label{G0 derivative est}
\frac{d}{dt}G_0(t,\vartheta)\Bigr|_{t=\pm\gamma}\ll\frac{\tan(2\vartheta)}{\sin^{3/2}(2\vartheta)}+\frac{1}{\sin^{1/2}(2\vartheta)}
\ll\frac{1}{\vartheta^{1/2}},
\end{equation}
which completes the proof of  \eqref{H0 estimate}.
\end{proof}
Substituting \eqref{H0 estimate} into \eqref{I Int4} we obtain for $\kc^{1/2-\delta}<r\le\kc/4$ (see \eqref{teta conditions}) that
\begin{equation}\label{I Int5}
I=\exp(i\kc \rho)\frac{2a_0\pi}{\kc^{1/3}}Ai(-\gamma^2\kc^{2/3})+
O\left(\frac{1}{(r\kc)^{1/2}}\right).
\end{equation}
Since $Ai(-x)\ll\min(1,x^{-1/4})$ we have
\begin{equation}\label{I Int6}
I\ll\frac{1}{\kc^{1/3}}\min\left(1,\frac{1}{\gamma^{1/2}\kc^{1/6}}\right)+
\frac{1}{(r\kc)^{1/2}}.
\end{equation}
It follows from \eqref{varteta def} and \eqref{gamma rho expression2}  that
\begin{multline*}
\gamma=2^{1/3}\left(\frac{\pi}{4}-\vartheta\right)+O\left((\pi/4-\vartheta)^5\right)=
2^{-2/3}\arccos\frac{4r}{\kc}+O\left(\arccos^5\frac{4r}{\kc}\right)=\\=
2^{-1/6}\left(1-\frac{4r}{\kc}\right)^{1/2}+O\left((1-4r/\kc)^{3/2}\right).
\end{multline*}
Consequently,
\begin{equation}\label{I Int7}
I\ll\frac{1}{\kc^{1/3}}\min\left(1,\frac{1}{(1-4r/\kc)^{1/4}\kc^{1/6}}\right)+
\frac{1}{(r\kc)^{1/2}}.
\end{equation}

Using \eqref{I Int7} to estimate the integral in \eqref{gMellin bettaInt}, we finally prove Lemma \ref{Lem gMellin small r}.



\section{Explicit formula for the mixed moment}\label{sec:5}

This section is devoted to proving an explicit formula for the mixed moment
\begin{equation*}
\sum_{\mathfrak{f}\in H_{4k}} \omega(\mathfrak{f}) L(\mathfrak{f},1/2)L(\sym^2\mathfrak{f}, 1/2).
\end{equation*}
To this end, we introduce two complex variables  $u$, $v$ with sufficiently large real parts
and consider the shifted moment
\begin{equation}
\M(u,v)=\sum_{\mathfrak{f}\in H_{4k}} \omega(\mathfrak{f}) L(\mathfrak{f},1/2+v)L(\sym^2\mathfrak{f}, 1/2+u).
\end{equation}
This  enables us to use the technique of analytic continuation.

Let us assume for simplicity that $0<\Re{u}<1$ and $\Re{v}>3/4+\Re{u}/2.$ Using \eqref{L def} and Lemma \ref{lem:EF} we obtain
\begin{equation}\label{M(u,v) 1}
\M(u,v)=\M^D(u,v)+\M^{ND}(u,v)+\ET_1(u,v)+\ET_2(u,v),
\end{equation}
where
\begin{equation}\label{MD MND def}
\M^D(u,v)=\sum_{l=1}^{\infty}\frac{M^{D}(u,l^2)}{l^{1+2v}},\quad
\M^{ND}(u,v)=\sum_{l=1}^{\infty}\frac{M^{ND}(u,l)}{l^{1/2+v}},
\end{equation}
\begin{equation}\label{ET1 ET2 def}
\ET_1(u,v)=\sum_{l=1}^{\infty}\frac{ET_1(u,l)}{l^{1/2+v}},\quad
\ET_2(u,v)=\sum_{l=1}^{\infty}\frac{ET_2(u,l)}{l^{1/2+v}}.
\end{equation}
As a consequence of \eqref{eq:MT} and \eqref{MD MND def}  we obtain
\begin{multline}\label{MD(u,v)}
\M^D(u,v)=
\zeta(3/2+2v+u)\zeta(1+2u) \\ + \zeta(3/2+2v-u)\zeta(1-2u)\sqrt{2}
(2\pi)^{3u}\cos{\pi(\frac{1}{4}+\frac{u}{2})}\\ \times
\frac{\Gamma(2k-1/4-u/2) \Gamma(2k+1/4-u/2)
\Gamma(1-2u)}{\Gamma(2k+1/4+u/2) \Gamma(2k-1/4+u/2)\Gamma(1-u)},
\end{multline}
\begin{multline}\label{MD(0,v)}
\M^D(0,v)=
\frac{\zeta(3/2+2v)}{2}
\Biggl(-3\log{2\pi}+\frac{\pi}{2}+3\gamma\\+
2\frac{\zeta'(3/2+2v)}{\zeta(3/2+2v)}
+\frac{\Gamma'}{\Gamma}(2k-1/4)+\frac{\Gamma'}{\Gamma}(2k+1/4)\Biggr).
\end{multline}
Similarly, it follows from \eqref{eq:lg}, \eqref{eq:MNDT}  and \eqref{MD MND def} that
\begin{equation}\label{MND(u,v)}
\M^{ND}(u,v)=
\frac{(2\pi)^{1/2+u}}{2}\frac{\Gamma(2k-1/4-u/2)}{\Gamma(2k+1/4+u/2)}
\sum_{l=1}^{\infty}\frac{\mathscr{L}_{-4l}(1/2+u)}{l^{3/4+v-u/2}},
\end{equation}
\begin{equation}\label{MND(0,v)}
\M^{ND}(0,v)=
\frac{2^{3/2}\pi}{\Gamma(3/4)}\frac{\Gamma(2k-1/4)}{\Gamma(2k+1/4)}L^{-}_{g}(1/4+v).
\end{equation}

We remark that a part of the main term is also contained in $\ET_1(u,v)$ and $\ET_2(u,v)$, which we analyse in detail in the next two subsections.
\subsection{Analysis of $\ET_2(u,v)$}
\begin{lem}
For $0<\Re{u}<1,$ $\Re{v}>3/4+\Re{u}/2$ we have
\begin{multline}\label{ET2 0}
\ET_2(u,v)=
(2\pi)^{1/2+u}2^{1+2v}\frac{1}{2\pi i}\\ \times \int_{(\sigma)}
\left(
\sum_{\substack{n=1\\ n\equiv0(2)}}^{\infty}\sum_{\substack{m=1\\ m\equiv0(4)}}^{\infty}+
\sum_{\substack{n=1\\ n\equiv1(2)}}^{\infty}\sum_{\substack{m=1\\ m\equiv1(4)}}^{\infty}
\right)
\frac{\mathscr{L}_{m}(1/2+u)}{m^{1/2+v+s}n^{1/2-u-2s}}\hat{f}_{2k}(u,v;s)ds,
\end{multline}
where $1/2-\Re{v}<\sigma<-1/4-\Re{u}/2.$
\end{lem}
\begin{proof}
Substituting \eqref{eq:ET2} into \eqref{ET1 ET2 def} we obtain
\begin{equation}\label{ET2 1}
\ET_2(u,v)=
\sum_{l=1}^{\infty}\frac{(2\pi)^{1/2+u}}{l^{1/2+v}}
\sum_{n>2\sqrt{l}} \frac{\mathscr{L}_{n^2-4l}(1/2+u)}{n^{1/2-u}}I_{2k}\left(u;\frac{n}{l^{1/2}}\right).
\end{equation}
It follows from \eqref{eq:integralI} that $I_{2k}\left(u;x\right)\sim x^{1-4k}$ as $x\rightarrow\infty.$ Thus
\begin{equation*}
\sum_{l=1}^{\infty}\frac{1}{l^{1/2+v}}
\sum_{n>2\sqrt{l}} \frac{\mathscr{L}_{n^2-4l}(1/2+u)}{n^{1/2-u}}I_{2k}\left(u;\frac{n}{l^{1/2}}\right)\ll
\sum_{l=1}^{\infty}\frac{l^{1/4+\Re{u}/2}}{l^{1/2+\Re{v}}}.
\end{equation*}
And we see that the double series on the right-hand side of \eqref{ET2 1} converges absolutely provided that $\Re{v}>3/4+\Re{u}/2.$ Changing the order of summation in \eqref{ET2 1} and making the change of variables $m=n^2-4l$, we obtain
\begin{equation*}
\ET_2(u,v)=(2\pi)^{1/2+u}
\sum_{n=1}^{\infty}
\sum_{\substack{0<m<n^2\\ m\equiv n^2(4)}}\frac{\mathscr{L}_{m}(1/2+u)2^{1+2v}}{(n^2-m)^{1/2+v}n^{1/2-u}}
I_{2k}\left(u;\frac{2n}{(n^2-m)^{1/2}}\right).
\end{equation*}
Rewriting this using \eqref{def f} yields
\begin{multline}\label{ET2 2}
\ET_2(u,v)=(2\pi)^{1/2+u}2^{1+2v}
\left(
\sum_{\substack{n=1\\ n\equiv0(2)}}^{\infty}\sum_{\substack{m=1\\ m\equiv0(4)}}^{\infty}+
\sum_{\substack{n=1\\ n\equiv1(2)}}^{\infty}\sum_{\substack{m=1\\ m\equiv1(4)}}^{\infty}
\right)\\\times
\frac{\mathscr{L}_{m}(1/2+u)}{m^{1/2+v}m^{1/2-u}}
f_{2k}\left(u,v;\frac{m}{n^2}\right).
\end{multline}
Applying the Mellin inversion formula for $f_{2k}\left(u,v;m/n^2\right)$ completes the proof.
\end{proof}

\begin{lem}
For $\Re{v}>3/4$ we have
\begin{equation}\label{ET2(0,v)}
\ET_2(0,v)=\frac{1}{2\pi i}\int_{(\sigma)}F_{2k}(v,s)ds,
\end{equation}
where $1/2-\Re{v}<\sigma<-1/4$  and
\begin{multline}\label{F def}
F_{2k}(v,s)=
(2\pi)^{1/2}2^{1+2v}
\Biggl(
\left(1-2^{2s-1/2}\right)\frac{2\sqrt{\pi}}{\Gamma(1/4)}L^{+}_{f}(s+v)-\\
-\left(1-2^{2s+1/2}\right)\frac{4\pi^{1/2}}{2^{1+2v+2s}\Gamma(1/4)}L^{+}_{g}(s+v)
\Biggr)
\zeta(1/2-2s)\hat{f}_{2k}(0,v;s).
\end{multline}
\end{lem}
\begin{proof}
We first let $u=0$ in \eqref{ET2 0}. Then for $\Re{v}>3/4$ the following formula holds
\begin{multline}\label{ET2(0,v)1}
\ET_2(0,v)=
(2\pi)^{1/2}2^{1+2v}\frac{1}{2\pi i}\int_{(\sigma)}
\left(
\sum_{\substack{n=1\\ n\equiv0(2)}}^{\infty}\sum_{\substack{m=1\\ m\equiv0(4)}}^{\infty}+
\sum_{\substack{n=1\\ n\equiv1(2)}}^{\infty}\sum_{\substack{m=1\\ m\equiv1(4)}}^{\infty}
\right)\\\times
\frac{\mathscr{L}_{m}(1/2)}{m^{1/2+v+s}n^{1/2-2s}}\hat{f}_{2k}(0,v;s)ds,
\end{multline}
where $1/2-\Re{v}<\sigma<-1/4.$
It follows from \eqref{eq:lg} that
\begin{equation}\label{Lg 1}
\sum_{\substack{m=1\\ m\equiv0(4)}}^{\infty}\frac{\mathscr{L}_{m}(1/2)}{m^{1/2+v+s}}=
\frac{4\pi^{1/2}}{4^{1/2+v+s}\Gamma(1/4)}L^{+}_{g}(s+v).
\end{equation}
Since $ \mathscr{L}_n(s) $ vanishes if $n \equiv 2,3 \Mod{4} $, we obtain using \eqref{eq:lf} and \eqref{Lg 1} that
\begin{multline}\label{Lg 2}
\sum_{\substack{m=1\\ m\equiv1(4)}}^{\infty}\frac{\mathscr{L}_{m}(1/2)}{m^{1/2+v+s}}=
\sum_{m=1}^{\infty}\frac{\mathscr{L}_{m}(1/2)}{m^{1/2+v+s}}-
\sum_{\substack{m=1\\ m\equiv0(4)}}^{\infty}\frac{\mathscr{L}_{m}(1/2)}{m^{1/2+v+s}}\\=
\frac{2\sqrt{\pi}}{\Gamma(1/4)}L^{+}_{f}(s+v)-\frac{4\pi^{1/2}}{4^{1/2+v+s}\Gamma(1/4)}L^{+}_{g}(s+v).
\end{multline}
Furthermore,
\begin{equation}\label{zeta21}
\sum_{\substack{n=1\\ n\equiv0(2)}}^{\infty}\frac{1}{n^{1/2-2s}}=\frac{\zeta(1/2-2s)}{2^{1/2-2s}},
\end{equation}
\begin{equation}\label{zeta22}
\sum_{\substack{n=1\\ n\equiv1(2)}}^{\infty}\frac{1}{n^{1/2-2s}}=\left(1-\frac{1}{2^{1/2-2s}}\right)\zeta(1/2-2s).
\end{equation}
Substituting \eqref{Lg 1}, \eqref{Lg 2}, \eqref{zeta21} and \eqref{zeta22} into \eqref{ET2(0,v)1} we prove \eqref{ET2(0,v)}.
\end{proof}
\begin{lem}
For $1/2<\Re{v}<3/4$ we have
\begin{multline}\label{ET2(0,v)1/2 3/4}
\ET_2(0,v)=
\frac{1}{2\pi i}\int_{(0)}F_{2k}(v,s)ds+\\ \sqrt{2\pi}2^{2v}\frac{\Gamma^2(1/4+v)}{\Gamma(1/4)}\frac{\Gamma(2k-1/4)\Gamma(2k-v)}{\Gamma(2k+1/4)\Gamma(2k+v)}
L^{+}_{f}(v-1/4),
\end{multline}
where $F_{2k}(v,s)$ is defined by \eqref{F def}.
\end{lem}
\begin{proof}
The function $F_{2k}(v,s)$ has a simple pole at $s=-1/4$ coming from $\zeta(1/2-2s).$ Moving the line of integration in \eqref{ET2(0,v)} to $\sigma=0$ we cross this pole, obtaining
\begin{equation}\label{ET2(0,v) 3}
\ET_2(0,v)=-\res_{s=-1/4}F_{2k}(v,s)+
\frac{1}{2\pi i}\int_{(0)}F_{2k}(v,s)ds.
\end{equation}
The right-hand side of \eqref{ET2(0,v) 3} shows that $\ET_2(0,v)$ can be continued to the region $\Re{v}>1/2.$  To prove \eqref{ET2(0,v)1/2 3/4} it remains to evaluate the residue. Using \eqref{F def} we have
\begin{equation*}
\res_{s=-1/4}F_{2k}(v,s)=\frac{2^{1/2+2v}\pi}{\Gamma(1/4)}\hat{f}_{2k}(0,v;-1/4)L^{+}_{f}(v-1/4).
\end{equation*}
Applying \eqref{fMellin -1/4} we prove the lemma.
\end{proof}

\begin{lem}
For $0\le\Re{v}<1/2$ we have
\begin{multline}\label{ET2(0,v)0 1/2}
\ET_2(0,v)=\res_{s=1/2-v}F_{2k}(v,s)+
\frac{1}{2\pi i}\int_{(0)}F_{2k}(v,s)ds\\
+\sqrt{2\pi}2^{2v}\frac{\Gamma^2(1/4+v)}{\Gamma(1/4)}\frac{\Gamma(2k-1/4)\Gamma(2k-v)}{\Gamma(2k+1/4)\Gamma(2k+v)}
L^{+}_{f}(v-1/4),
\end{multline}
where $F_{2k}(v,s)$ is defined by \eqref{F def}.
\end{lem}
\begin{proof}
The function $F_{2k}(v,s)$ has a double pole at $s=1/2-v$ coming from $L^{+}_{f,g}(s+v).$  To prove the analytic continuation of $\ET_2(0,v)$ to the region $\Re{v}<1/2$, we apply \cite[Corollary 2.4.2, p. 55]{CMR}, which yields for $\Re{v}<1/2$ that
\begin{multline*}
\ET_2(0,v)=-\res_{s=-1/4}F_{2k}(v,s)+\res_{s=1/2-v}F_{2k}(v,s)+
\frac{1}{2\pi i}\int_{(0)}F_{2k}(v,s)ds.
\end{multline*}
\end{proof}

\begin{lem}\label{lem:second type}
The following formula holds
\begin{equation}\label{ET2(0,0)expression1}
\ET_2(0,0)=\M^{ND}(0,0)+\res_{s=1/2-v}F_{2k}(v,s)+
\frac{1}{2\pi i}\int_{(0)}F_{2k}(v,s)ds.
\end{equation}
\end{lem}
\begin{proof}
Comparing \eqref{MND(0,v)} and \eqref{ET2(0,v)0 1/2}, we find that in order to
prove \eqref{ET2(0,0)expression1}, it is required to show that
\begin{equation*}
\sqrt{2\pi}\Gamma(1/4)L^{+}_{f}(-1/4)=
\frac{2^{3/2}\pi}{\Gamma(3/4)}L^{-}_{g}(1/4).
\end{equation*}
 Since $\Gamma(1/4)\Gamma(3/4)=\pi\sqrt{2}$ we need to verify that
\begin{equation*}
\sqrt{\pi}L^{+}_{f}(-1/4)=2^{1/2}L^{-}_{g}(1/4),
\end{equation*}
and this follows from \eqref{eq:lgs}.
\end{proof}



\begin{lem}
For any $\epsilon>0$ we have
\begin{equation}\label{ET2(0,0)expression2}
\ET_2(0,0)=\M^{ND}(0,0)+O\left(\frac{k^{\epsilon}}{k}\right),
\end{equation}
\begin{equation}\label{ET2(0,0) estimate}
\ET_2(0,0)\ll k^{-1/2}.
\end{equation}
\end{lem}
\begin{proof}
This follows immediately from \eqref{ET2(0,v)0 1/2}, \eqref{F def},  Lemma \ref{Lem fMellin 1/2},
Theorem \ref{thm second mom Lfg} and the estimate \eqref{fMellin ir}.
\end{proof}

\subsection{Analysis of $\ET_1(u,v)$}
\begin{lem}
Assume that  $0<\Re{u}<1$ and  $\Re{v}>3/4+\Re{u}/2+\max(\theta(1-2\Re{u}),0)$. Then the following formula holds
\begin{multline}\label{ET1 0}
\ET_1(u,v)=
(2\pi)^{1/2+u}2^{1+2v} \frac{1}{2\pi i}\times \\ \int_{(\sigma)}
\left(
\sum_{\substack{n=1\\ n\equiv0(2)}}^{\infty}\sum_{\substack{m=1\\ m\equiv0(4)}}^{\infty}+
\sum_{\substack{n=1\\ n\equiv1(2)}}^{\infty}\sum_{\substack{m=-1\\ m\equiv1(4)}}^{\infty}
\right)
\frac{\mathscr{L}_{-m}(1/2+u)}{m^{1/2+v+s}n^{1/2-u-2s}}\hat{g}_{2k}(u,v;s)ds,
\end{multline}
where $1/2-\Re{v}<\sigma<-1/4-\Re{u}/2.$
\end{lem}
\begin{proof}
Substituting \eqref{eq:ET1} into \eqref{ET1 ET2 def} we show that
\begin{equation}\label{ET1 1}
\ET_1(u,v)=
\sum_{l=1}^{\infty}\frac{(2\pi)^{1/2+u}}{l^{1/2+v}}
\sum_{0<n<2\sqrt{l}}\frac{\mathscr{L}_{n^2-4l}(1/2+u)}{n^{1/2-u}}I_{2k}\left(u;\frac{n}{l^{1/2}}\right).
\end{equation}
It follows from \eqref{eq:integralI} that $I_{2k}\left(u;x\right)\sim x^{1/2-\Re{u}}$ as $x\rightarrow0.$ Using this fact and applying \eqref{eq:subcL}, we obtain
\begin{multline*}
\sum_{l=1}^{\infty}\frac{1}{l^{1/2+v}}
\sum_{0<n<2\sqrt{l}}\frac{\mathscr{L}_{n^2-4l}(1/2+u)}{n^{1/2-u}}I_{2k}\left(u;\frac{n}{l^{1/2}}\right)\ll
\sum_{l=1}^{\infty}\frac{1}{l^{1/2+\Re{v}}}\frac{l^{1/2+\max(\theta(1-2\Re{u}),0)+\epsilon}}{l^{1/4-\Re{u}/2}}.
\end{multline*}
Therefore, the double series on the right-hand side of \eqref{ET1 1} converges absolutely provided that $\Re{v}>3/4+\Re{u}/2+\max(\theta(1-2\Re{u}),0)$.
Changing the order of summation in \eqref{ET1 1} and making the change of variables $-m=n^2-4l$, we have
\begin{equation*}
\ET_1(u,v)=(2\pi)^{1/2+u}
\sum_{n=1}^{\infty}
\sum_{\substack{m=1\\ m+n^2\equiv0(4)}}\frac{\mathscr{L}_{-m}(1/2+u)2^{1+2v}}{(n^2+m)^{1/2+v}n^{1/2-u}}
I_{2k}\left(u;\frac{2n}{(n^2+m)^{1/2}}\right).
\end{equation*}
Applying \eqref{def g}, we obtain
\begin{multline}\label{ET1 2}
\ET_1(u,v)=(2\pi)^{1/2+u}2^{1+2v}
\left(
\sum_{\substack{n=1\\ n\equiv0(2)}}^{\infty}\sum_{\substack{m=1\\ m\equiv0(4)}}^{\infty}+
\sum_{\substack{n=1\\ n\equiv1(2)}}^{\infty}\sum_{\substack{m=1\\ m\equiv-1(4)}}^{\infty}
\right)\\\times
\frac{\mathscr{L}_{-m}(1/2+u)}{m^{1/2+v}m^{1/2-u}}
g_{2k}\left(u,v;\frac{m}{n^2}\right).
\end{multline}
Using the Mellin inversion formula for $g_{2k}\left(u,v;m/n^2\right)$ we prove the lemma.
\end{proof}


\begin{lem}
For $\Re{v}>3/4$ the following representation takes place
\begin{equation}\label{ET1(0,v)3/4}
\ET_1(0,v)=\frac{1}{2\pi i}\int_{(\sigma)}G_{2k}(v,s)ds,
\end{equation}
where $1/2-\Re{v}<\sigma<-1/4$ and
\begin{multline}\label{G def}
G_{2k}(v,s)=
(2\pi)^{1/2}2^{1+2v}
\Biggl(
\left(1-2^{2s-1/2}\right)\frac{2\sqrt{\pi}}{\Gamma(3/4)}L^{-}_{f}(s+v)-\\-
\left(1-2^{2s+1/2}\right)\frac{4\pi^{1/2}}{2^{1+2v+2s}\Gamma(3/4)}L^{-}_{g}(s+v)
\Biggr)
\zeta(1/2-2s)\hat{g}_{2k}(0,v;s).
\end{multline}
\end{lem}
\begin{proof}
Letting $u=0$ in \eqref{ET1 0}, we obtain for $\Re{v}>3/4+\theta$
\begin{multline}\label{ET1(0,v)1}
\ET_1(0,v)=
(2\pi)^{1/2}2^{1+2v}\frac{1}{2\pi i}\times \\
\int_{(\sigma)}
\left(
\sum_{\substack{n=1\\ n\equiv0(2)}}^{\infty}\sum_{\substack{m=1\\ m\equiv0(4)}}^{\infty}+
\sum_{\substack{n=1\\ n\equiv1(2)}}^{\infty}\sum_{\substack{m=1\\ m\equiv-1(4)}}^{\infty}
\right)
\frac{\mathscr{L}_{-m}(1/2)}{m^{1/2+v+s}n^{1/2-2s}}\hat{g}_{2k}(0,v;s)ds,
\end{multline}
where $1/2-\Re{v}<\sigma<-1/4.$
It follows from \eqref{eq:lg} that
\begin{equation}\label{Lg 3}
\sum_{\substack{m=1\\ m\equiv0(4)}}^{\infty}\frac{\mathscr{L}_{-m}(1/2)}{m^{1/2+v+s}}=
\frac{4\pi^{1/2}}{4^{1/2+v+s}\Gamma(3/4)}L^{-}_{g}(s+v).
\end{equation}
Recall that $ \mathscr{L}_n(s) $ vanishes for $n \equiv 2,3 \Mod{4} $. Consequently, using \eqref{eq:lf} and \eqref{Lg 3} we show that
\begin{multline}\label{Lg 4}
\sum_{\substack{m=1\\ m\equiv-1(4)}}^{\infty}\frac{\mathscr{L}_{-m}(1/2)}{m^{1/2+v+s}}=
\sum_{m=1}^{\infty}\frac{\mathscr{L}_{-m}(1/2)}{m^{1/2+v+s}}-
\sum_{\substack{m=1\\ m\equiv0(4)}}^{\infty}\frac{\mathscr{L}_{-m}(1/2)}{m^{1/2+v+s}}\\=
\frac{2\sqrt{\pi}}{\Gamma(3/4)}L^{-}_{f}(s+v)-\frac{4\pi^{1/2}}{4^{1/2+v+s}\Gamma(3/4)}L^{-}_{g}(s+v).
\end{multline}
Applying \eqref{zeta21}, \eqref{zeta22}, \eqref{Lg 3} and \eqref{Lg 4}  to evaluate \eqref{ET1(0,v)1}, we prove \eqref{ET1(0,v)3/4}.
\end{proof}


\begin{lem}
For $1/2<\Re{v}<3/4$ the following formula holds
\begin{multline}\label{ET1(0,v)1/2 3/4}
\ET_1(0,v)=
\frac{1}{2\pi i}\int_{(0)}G_{2k}(v,s)ds
\\-\sqrt{\pi}2^{1+2v}\sin(\pi v)
\frac{\Gamma^2(1/4+v)}{\Gamma(3/4)}\frac{\Gamma(2k-1/4)\Gamma(2k-v)}{\Gamma(2k+1/4)\Gamma(2k+v)}
L^{-}_{f}(v-1/4),
\end{multline}
where $G_{2k}(v,s)$ is defined by \eqref{G def}.
\end{lem}
\begin{proof}
The function $G_{2k}(v,s)$ has a simple pole at $s=-1/4$ from $\zeta(1/2-2s).$  Moving the line of integration in \eqref{ET1(0,v)3/4} to $\sigma=0$ we cross this pole, getting
\begin{equation}\label{ET1(0,v) 3}
\ET_1(0,v)=-\res_{s=-1/4}G_{2k}(v,s)+
\frac{1}{2\pi i}\int_{(0)}G_{2k}(v,s)ds.
\end{equation}
The right-hand side of \eqref{ET1(0,v) 3} proves the analytic continuation of $\ET_1(0,v)$ to the region $\Re{v}>1/2.$
Then to complete the proof of \eqref{ET1(0,v)1/2 3/4}, it remains to evaluate the residue. Using \eqref{G def} we have
\begin{equation*}
\res_{s=-1/4}G_{2k}(v,s)=\frac{2^{1/2+2v}\pi}{\Gamma(1/4)}\hat{g}_{2k}(0,v;-1/4)L^{-}_{f}(v-1/4).
\end{equation*}
The lemma follows by applying \eqref{gMellin -1/4}.
\end{proof}

\begin{lem}\label{lem:et1at0v}
For $0\le\Re{v}<1/2$ we have 
\begin{multline}\label{ET1(0,v)0 1/2}
\ET_1(0,v)=\res_{s=1/2-v}G_{2k}(v,s)+\frac{1}{2\pi i}\int_{(0)}G_{2k}(v,s)ds
\\-\sqrt{\pi}2^{1+2v}\sin(\pi v)
\frac{\Gamma^2(1/4+v)}{\Gamma(1/4)}\frac{\Gamma(2k-1/4)\Gamma(2k-v)}{\Gamma(2k+1/4)\Gamma(2k+v)}
L^{+}_{f}(v-1/4),
\end{multline}
where $G_{2k}(v,s)$ is defined by \eqref{G def} and 
\begin{equation}\label{res G 1/2-v}
\res_{s=1/2-v}G_{2k}(v,s)=M_1(v)+M_2(v),
\end{equation}
\begin{multline}\label{res G 1/2-v M1}
M_1(v)=
\frac{2^{5/2+2v}\pi}{\Gamma(3/4)}
\Biggl[
c_f^{-}(-1)\zeta(2v-1/2)\hat{g}_{2k}(0,v;1/2-v)\left(1-2^{1/2-2v}\right)\\+
c_f^{-}(-2)
\Biggl(
\zeta(2v-1/2)\left(1-2^{1/2-2v}\right)\frac{\partial}{\partial s}\hat{g}_{2k}(0,v;s)\Bigg|_{s=1/2-v}\\-
2\zeta'(2v-1/2)\hat{g}_{2k}(0,v;1/2-v)\left(1-2^{1/2-2v}\right)\\ -
\zeta(2v-1/2)\hat{g}_{2k}(0,v;1/2-v)2^{3/2-2v}\log2
\Biggr)\Biggr],
\end{multline}
\begin{multline}\label{res G 1/2-v M2}
M_2(v)=
\frac{8\pi}{\Gamma(3/4)}
\Biggl[
c_g^{-}(-1)\zeta(2v-1/2)\hat{g}_{2k}(0,v;1/2-v)\left(1-2^{2v-3/2}\right)\\+
c_g^{-}(-2)
\Biggl(
\zeta(2v-1/2)\left(1-2^{2v-3/2}\right)\frac{\partial}{\partial s}\hat{g}_{2k}(0,v;s)\Bigg|_{s=1/2-v}\\-
2\zeta'(2v-1/2)\hat{g}_{2k}(0,v;1/2-v)\left(1-2^{2v-3/2}\right)\\+
\zeta(2v-1/2)\hat{g}_{2k}(0,v;1/2-v)2^{2v-1/2}\log2
\Biggr)
\Biggr].
\end{multline}
\end{lem}
\begin{proof}
The function $G_{2k}(v,s)$ has a double pole at $s=1/2-v$ from $L^{-}_{f,g}(s+v).$
To prove  the analytic continuation of $\ET_1(0,v)$ to the region $\Re{v}<1/2$ we apply \cite[Corollary 2.4.2, p. 55]{CMR}. Consequently, for $\Re{v}<1/2$ we have
\begin{multline}\label{ET1(0,v) 4}
\ET_1(0,v)=\res_{s=1/2-v}G_{2k}(v,s)-\res_{s=-1/4}G_{2k}(v,s)+\\
\frac{1}{2\pi i}\int_{(0)}G_{2k}(v,s)ds.
\end{multline}
Then it follows from  \eqref{G def} that
\begin{multline}\label{res G 1/2-v 0}
\res_{s=1/2-v}G_{2k}(v,s)=\\
\frac{2^{5/2+2v}\pi}{\Gamma(3/4)}
\res_{s=1/2-v}\Biggl(
\zeta(1/2-2s)\hat{g}_{2k}(0,v;s)\left(1-2^{2s-1/2}\right)L^{-}_{f}(s+v)
\Biggr)\\+
\frac{8\pi}{\Gamma(3/4)}
\res_{s=1/2-v}\Biggl(
\zeta(1/2-2s)\hat{g}_{2k}(0,v;s)\left(1-2^{-2s-1/2}\right)L^{-}_{g}(s+v)
\Biggr).
\end{multline}
Let $H(s)$ be an arbitrary function that is holomorphic at $s=1/2-v$. Using \eqref{Laurent series Lfg}, we obtain the Laurent series
\begin{multline}\label{Laurent series1}
H(s)L^{-}_{f,g}(s+v)=\frac{c_{f,g}^{-}(-2)H(1/2-v)}{(s+v-1/2)^2}\\+\frac{c_{f,g}^{-}(-1)H(1/2-v)+c_{f,g}^{-}(-2)H'(1/2-v)}{s+v-1/2}+O(1).
\end{multline}
Applying \eqref{Laurent series1} to evaluate \eqref{res G 1/2-v 0}, we prove \eqref{res G 1/2-v}.
\end{proof}
\subsection{Analytic continuation}

Finally, we obtain the following decomposition for the mixed moment.
\begin{thm}
For $\Re{v}\ge0$ we have
\begin{equation}\label{M(u,v) exact formula}
\M(0,v)=\M^D(0,v)+\M^{ND}(0,v)+\ET_1(0,v)+\ET_2(0,v),
\end{equation}
where $\M^D(0,v)$ is defined by \eqref{MD(0,v)} and $\M^{ND}(0,v)$ by \eqref{MND(0,v)}. Furthermore, the terms $\ET_1(0,v)$ and  $\ET_2(0,v)$ are given by
\eqref{ET1(0,v)3/4} and  \eqref{ET2(0,v)} for $\Re{v}>3/4$, by \eqref{ET1(0,v)1/2 3/4} and \eqref{ET2(0,v)1/2 3/4} for
$1/2<\Re{v}\le3/4$ and by \eqref{ET1(0,v)0 1/2} and \eqref{ET2(0,v)0 1/2} for $0\le\Re{v}<1/2.$
\end{thm}
\begin{proof}
In order to prove the theorem, it remains to show that the right-hand side of \eqref{M(u,v) exact formula} is holomorphic for $\Re{v}\ge0.$  More precisely, we need to consider points $v=3/4$ and $v=1/4.$ The only summands on the right-hand side of \eqref{M(u,v) exact formula} that are not holomorphic at $v=3/4$ come from \eqref{ET1(0,v)1/2 3/4} and \eqref{ET2(0,v)1/2 3/4}, namely:
\begin{multline}\label{res at 3/4}
\sqrt{2\pi}2^{2v}\frac{\Gamma^2(1/4+v)}{\Gamma(1/4)}\frac{\Gamma(2k-1/4)\Gamma(2k-v)}{\Gamma(2k+1/4)\Gamma(2k+v)}
L^{+}_{f}(v-1/4)-\\
-\sqrt{\pi}2^{1+2v}\sin(\pi v)
\frac{\Gamma^2(1/4+v)}{\Gamma(3/4)}\frac{\Gamma(2k-1/4)\Gamma(2k-v)}{\Gamma(2k+1/4)\Gamma(2k+v)}
L^{-}_{f}(v-1/4)=\\=
\sqrt{\pi}2^{2v}\Gamma^2(1/4+v)\frac{\Gamma(2k-1/4)\Gamma(2k-v)}{\Gamma(2k+1/4)\Gamma(2k+v)}\times \\
\left(\frac{\sqrt{2}}{\Gamma(1/4)}L^{+}_{f}(v-1/4)-
\frac{2\sin(\pi v)}{\Gamma(3/4)}L^{-}_{f}(v-1/4)
\right).
\end{multline}
Therefore, to prove that the right-hand side of \eqref{M(u,v) exact formula} is holomorphic at $v=3/4$, it is sufficient to show that
\begin{equation}
\frac{\sqrt{2}}{\Gamma(1/4)}L^{+}_{f}(v-1/4)-
\frac{2\sin(\pi v)}{\Gamma(3/4)}L^{-}_{f}(v-1/4)
\end{equation}
is holomorphic at $v=3/4.$ Using Theorem \ref{cor Lfg coeff relation} and the asymptotic formula
$$\sin(\pi v)=\frac{1}{\sqrt{2}}-\frac{\pi(v-3/4)}{\sqrt{2}}+O((v-3/4)^2),$$
we obtain
\begin{multline}
\frac{\sqrt{2}}{\Gamma(1/4)}L^{+}_{f}(v-1/4)-
\frac{2\sin(\pi v)}{\Gamma(3/4)}L^{-}_{f}(v-1/4)=\\=
\frac{1}{(v-3/4)^2}\left(
\frac{c^{+}_{f}(-2)\sqrt{2}}{\Gamma(1/4)}-
\frac{c^{-}_{f}(-2)\sqrt{2}}{\Gamma(3/4)}
\right)+\\+
\frac{1}{v-3/4}\left(
\frac{c^{+}_{f}(-1)\sqrt{2}}{\Gamma(1/4)}-
\frac{c^{-}_{f}(-1)\sqrt{2}}{\Gamma(3/4)}+
\frac{c^{-}_{f}(-2)\pi\sqrt{2}}{\Gamma(3/4)}
\right)+O(1)=O(1).
\end{multline}
 Thus the right-hand side of \eqref{M(u,v) exact formula} is holomorphic at $v=3/4$.

The only summands on  the right-hand side of \eqref{M(u,v) exact formula} that are not holomorphic at $v=1/4$ come from \eqref{MND(0,v)} and  \eqref{res G 1/2-v}, namely
\begin{equation*}
\M^{ND}(0,v)+\res_{s=1/2-v}G_{2k}(v,s).
\end{equation*}

Let us consider the $M_1(v)$ part of $\res_{s=1/2-v}G_{2k}(v,s)$ given by \eqref{res G 1/2-v M1}. Using \eqref{gMellin v 1/2-v} and \eqref{gMellin dif v 1/2-v}, we obtain a Laurent series for $M_1(v)$  at the point $v=1/4$:
\begin{multline}\label{M1 1/4 1}
M_1(v)=
\frac{2^{5/2}\pi}{\Gamma(3/4)}
c_f^{-}(-2)
\Biggl(
\zeta(2v-1/2)\left(2^{2v}-2^{1/2}\right)\frac{\partial}{\partial s}\hat{g}_{2k}(0,v;s)\Bigg|_{s=1/2-v}\\-
\frac{2^{3}\zeta(0)\log2}{2v-1/2}\frac{\Gamma(2k-1/4)}{\Gamma(2k+1/4)}
\Biggr)+O(1).
\end{multline}
Using \eqref{gMellin dif v 1/2-v} and the fact that
\begin{equation*}
\zeta(2v-1/2)\left(2^{2v}-2^{1/2}\right)=(2v-1/2)\zeta(0)2^{1/2}\log2+O((2v-1/2)^2),
\end{equation*}
  we show that $M_1(v)=O(1)$ as $v\to 1/4.$

Let us consider the $M_2(v)$ part of $\res_{s=1/2-v}G_{2k}(v,s)$ given by \eqref{res G 1/2-v M2}. Applying \eqref{gMellin v 1/2-v},  we obtain
\begin{multline}\label{M2 Laurent 1/4 1}
M_2(v)=\frac{8\pi}{\Gamma(3/4)}c_g^{-}(-2)
\zeta(2v-1/2)\left(1-2^{2v-3/2}\right)\frac{\partial}{\partial s}\hat{g}_{2k}(0,v;s)\Bigg|_{s=1/2-v}\\+
\frac{1}{2v-1/2}\frac{8\pi}{\Gamma(3/4)}\frac{\Gamma(2k-1/4)}{\Gamma(2k+1/4)}
\Biggl(
c_g^{-}(-1)\zeta(0)2^{1/2}\\+
c_g^{-}(-2)
\Biggl(
-2^{3/2}\zeta'(0)+
\zeta(0)2^{3/2}\log2
\Biggr)
\Biggr).
\end{multline}
In order to evaluate a Laurent series for the remaining term we use \eqref{gMellin dif v 1/2-v}  together with the following formula
\begin{equation*}
\zeta(2v-1/2)\left(1-2^{2v-3/2}\right)=\frac{\zeta(0)}{2}+
\frac{\zeta'(0)-\zeta(0)\log2}{2}(2v-1/2)+O((2v-1/2)^2).
\end{equation*}
Consequently,
\begin{multline}\label{M2 Laurent 1/4 2}
M_2(v)=\frac{8\pi}{\Gamma(3/4)}\frac{\Gamma(2k-1/4)}{\Gamma(2k+1/4)}
\frac{c_g^{-}(-2)2^{3/2}\zeta(0)}{(2v-1/2)^2}
+\\+
\frac{8\pi}{\Gamma(3/4)}\frac{\Gamma(2k-1/4)}{\Gamma(2k+1/4)}
\frac{c_g^{-}(-1)2^{1/2}\zeta(0)}{2v-1/2}+O(1).
\end{multline}

It follows from \eqref{MND(0,v)} and \eqref{Laurent series Lfg} that
\begin{multline}\label{MND Laurent 1/4}
\M^{ND}(0,v)=\frac{2^{7/2}\pi}{\Gamma(3/4)}\frac{\Gamma(2k-1/4)}{\Gamma(2k+1/4)}
\frac{c_g^{-}(-2)}{(2v-1/2)^2}
+\\+
\frac{2^{5/2}\pi}{\Gamma(3/4)}\frac{\Gamma(2k-1/4)}{\Gamma(2k+1/4)}
\frac{c_g^{-}(-1)}{2v-1/2}
+O(1).
\end{multline}
Since $M_1(v)=O(1),$ applying \eqref{M2 Laurent 1/4 2} and \eqref{MND Laurent 1/4}
we conclude that the sum
\begin{equation*}
\M^{ND}(0,v)+\res_{s=1/2-v}G_{2k}(v,s),
\end{equation*} 
and consequently the right-hand side of \eqref{M(u,v) exact formula}, are holomorphic at $v=1/4.$ 
\end{proof}

\begin{lem}
The following asymptotic formula holds
\begin{equation}\label{resG=MD}
\res_{s=1/2}G_{2k}(0,s)=\M^{D}(0,0)+O(k^{-2}).
\end{equation} 
\end{lem}
\begin{proof}
We compare the leading terms. It follows from \cite[(5.11.2)]{HMF} that 
\begin{equation}\label{psi relation}
2\frac{\Gamma'}{\Gamma}(2k)=\frac{\Gamma'}{\Gamma}(2k-1/4)+\frac{\Gamma'}{\Gamma}(2k+1/4)+O(k^{-2}).
\end{equation} 
Therefore, \eqref{MD(0,v)} implies that the leading term of $\M^{D}(0,0)$ is equal to 
\begin{equation}\label{MD(0,0)leading term}
\zeta(3/2)\frac{\Gamma'}{\Gamma}(2k).
\end{equation}
 Let us compute the leading term of $\res_{s=1/2}G_{2k}(0,s).$ It follows from \eqref{gMellin v 1/2-v} and \eqref{gMellin dif v 1/2-v} that the leading term is 
\begin{multline}\label{resG leading term1}
2\frac{\Gamma'}{\Gamma}(2k)\frac{16\pi}{\Gamma(3/4)}\Biggl(
c_f^{-}(-2)(1-\sqrt{2})\zeta(-1/2)\Gamma(-1/2)+\\+
c_g^{-}(-2)(\sqrt{2}-1/2)\zeta(-1/2)\Gamma(-1/2)
\Biggr).
\end{multline}
Using the functional equation for the Riemann zeta function \cite[(25.4.1)]{HMF} we have $$\zeta(-1/2)\Gamma(-1/2)=\frac{\zeta(3/2)}{2\sqrt{\pi}}.$$
Consequently, the leading term of $\res_{s=1/2}G_{2k}(0,s)$ is as follows:
\begin{equation}\label{resG leading term2}
\frac{\Gamma'}{\Gamma}(2k)\zeta(3/2)\frac{16\sqrt{\pi}}{\Gamma(3/4)}\left(
c_f^{-}(-2)(1-\sqrt{2})+
c_g^{-}(-2)(\sqrt{2}-1/2)
\right).
\end{equation}
Finally, applying \eqref{relation for MT 1}, we find that  \eqref{resG leading term2} is equal to \eqref{MD(0,0)leading term}.
\end{proof}


\subsection{Proof of main theorems}
\begin{proof}[Proof of Theorem \ref{main thm}]
Asymptotic formula \eqref{M(0,0) exact formula1} is a direct consequence of \eqref{M(u,v) exact formula} for $v=0$. More precisely, we replace $\ET_2(0,v)$ by \eqref{ET2(0,0)expression2},
$\ET_1(0,v)$ by \eqref{ET1(0,v)0 1/2}, and apply \eqref{resG=MD}.
\end{proof}


\begin{proof}[Proof of Theorem \ref{maincor}]
Consider \eqref{M(0,0) exact formula1} and note that all summands except the integral can be trivially bounded by $\log k$. The final step is to show that
\begin{equation}\label{integral of G estimate}
\frac{1}{2\pi i}\int_{(0)}G_{2k}(0,s)ds\ll\log^{3} k,
\end{equation}
where $G_{2k}(0,s)$ is given by \eqref{G def}. In view of \eqref{G def}, in order to establish \eqref{integral of G estimate}, we need to prove that
\begin{equation}\label{integral of G estimate2}
\frac{1}{2\pi i}\int_{(0)}L^{-}_{f,g}(s)\zeta(1/2-2s)\hat{g}_{2k}(0,0;s)ds\ll\log^{3} k.
\end{equation}
Let $r=\Im{s}$. By Lemma \ref{Lem gMellin big r} the contribution of $|r|>3k$ is negligible. As a result, we are left to show that
\begin{equation}\label{integral of G estimate3}
\int_{-3k}^{3k}L^{-}_{f,g}(ir)\zeta(1/2-2ir)\hat{g}_{2k}(0,0;ir)dr\ll\log^{3} k.
\end{equation}
Let $\delta$ be some fixed constant such that $0<\delta<1/4$. For $|r|<\kc^{1/2-\delta}$ we apply the trivial bound \eqref{gMellin triv.est.}. Then it is required to establish that
\begin{equation}\label{integral of G estimate4}
\int_{\kc^{1/2-\delta}}^{3k}L^{-}_{f,g}(ir)\zeta(1/2-2ir)\hat{g}_{2k}(0,0;ir)dr\ll\log^{3} k.
\end{equation}
With this goal, we apply Lemma \ref{Lem gMellin small r}, Theorem \ref{thm second mom Lfg}, the following estimate for the second moment of the Riemann zeta function
\begin{equation}\label{zeta 2moment}
\int_{T}^{T+H}|\zeta(1/2+ir)|^2dr\ll H\log T
\end{equation}
over short intervals $H\gg T^{1/3}$, and  the Cauchy-Schwarz inequality. Consequently, we show that the contribution of the second summand on the right-hand side of
\eqref{gMellin small r est}  to \eqref{integral of G estimate4} is negligibly small. So our problem got reduced to proving that
\begin{equation}\label{integral of G estimate5}
\int_{\kc^{1/2-\delta}}^{3k}|L^{-}_{f,g}(ir)\zeta(1/2-2ir)|\frac{1}{\kc^{5/6}}\min\left(1,\frac{\kc^{1/12}}{|\kc-4r|^{1/4}}\right)dr\ll\log^{3} k.
\end{equation}
Opening the minimum we obtain three integrals:
\begin{multline}\label{integral of G estimate6}
\int_{r_1}^{r_2}\frac{|L^{-}_{f,g}(ir)\zeta(1/2-2ir)|}{\kc^{3/4}|\kc-4r|^{1/4}}dr+
\int_{r_2}^{r_3}|L^{-}_{f,g}(ir)\zeta(1/2-2ir)|\frac{dr}{\kc^{5/6}}+\\
\int_{r_3}^{r_4}\frac{|L^{-}_{f,g}(ir)\zeta(1/2-2ir)|}{\kc^{3/4}|\kc-4r|^{1/4}}dr
\ll\log^{3} k,
\end{multline}
where $r_1=\kc^{1/2-\delta}$, $r_2=\kc/4-\kc^{1/3}$, $r_3=\kc/4+\kc^{1/3},r_4=3k.$ To estimate the second integral we apply the Cauchy-Schwarz inequality, Theorem \ref{thm second mom Lfg} and \eqref{zeta 2moment}, getting 
\begin{equation}\label{integral of G estimate7}
\int_{r_2}^{r_3}|L^{-}_{f,g}(ir)\zeta(1/2-2ir)|\frac{dr}{\kc^{5/6}}
\ll\frac{\log^{5/2} k}{k^{1/6}}.
\end{equation}
Let us now consider the first integral in \eqref{integral of G estimate6}. Applying the Cauchy-Schwarz inequality and Theorem \ref{thm second mom Lfg}, we obtain
\begin{equation}\label{integral of G estimate8}
\int_{r_1}^{r_2}\frac{|L^{-}_{f,g}(ir)\zeta(1/2-2ir)|}{\kc^{3/4}|\kc-4r|^{1/4}}dr
\ll\frac{\log^{2} k}{\kc^{1/4}}\left(\int_{r_1}^{r_2}\frac{|\zeta(1/2-2ir)|^2}{|\kc-4r|^{1/2}}dr\right)^{1/2}.
\end{equation}
Making the change of variable $\kc-4r=x$, and then performing a dyadic partition of unity, we prove using \eqref{zeta 2moment} that
\begin{equation}\label{integral of G estimate9}
\int_{r_1}^{r_2}\frac{|\zeta(1/2-2ir)|^2}{|\kc-4r|^{1/2}}dr\ll k^{1/2}\log^2k,
\end{equation}
where one of the logarithms comes from the partition of unity. Substituting  \eqref{integral of G estimate9} into \eqref{integral of G estimate8} we obtain 
\begin{equation}\label{integral of G estimate10}
\int_{r_1}^{r_2}\frac{|L^{-}_{f,g}(ir)\zeta(1/2-2ir)|}{\kc^{3/4}|\kc-4r|^{1/4}}dr
\ll\log^{3} k.
\end{equation}
Finally, the third integral in \eqref{integral of G estimate6} can be estimated in the same way as the first one. This completes the proof.
\end{proof}




\end{document}